\documentclass[a4paper,12pt]{article}

\usepackage{latexsym}
\usepackage{amsthm}
\usepackage{amsmath}
\usepackage{amsfonts, epsfig, amsmath, amssymb, color}
\usepackage[cmtip,arrow]{xy}
\usepackage{pb-diagram, pb-xy}
\usepackage{amssymb,epsfig}

\newfont{\sdbl}{msbm9}
\newfont{\dbl}{msbm10 at 12pt}
\theoremstyle{definition}
\newcommand{\cl}{{{\cal L}}}

\newcommand{\cc}{{{\cal C}}}
\newcommand{\cn}{{{\cal N}}}
\newcommand{\ck}{{{\cal K}}}

\newcommand{\cf}{{{\cal F}}}
\newcommand{\cm}{{{\cal M}}}

\newcommand{\co}{{{\cal O}}}
\newcommand{\ca}{{{\cal A}}}
\newcommand{\cb}{{{\cal B}}}
\newcommand{\cg}{{{\cal G}}}
\newcommand{\h}{{{\cal H}}}
\newcommand{\G}{{{\bf G}}}
\newcommand{\xo}{{\mbox{\em \r{X}}}}
\newcommand{\sxo}{{\scriptsize  \mbox{\em \r{X}}}}
\newcommand{\ssxo}{{\tiny  \mbox{\em \r{X}}}}
\newcommand{\eo}{{\mbox{\em \r{E}}}}

\newcommand{\dz}{{\mbox{\dbl Z}}}

\newcommand{\dn}{{\mbox{\dbl N}}}

\newcommand{\dv}{{\mbox{\dbl V}}}

\newcommand{\sdz}{{\mbox{\sdbl Z}}}
\newcommand{\sdn}{{\mbox{\sdbl N}}}

\newcommand{\chara}{\mathop{\rm char}\nolimits}

\newcommand{\ord}{\mathop{\rm \bf ord}\nolimits}

\newcommand{\Ob}{\mathop{\rm Ob}\nolimits}
\newcommand{\bc}{{{\bf C}}}
\newcommand{\nn}{{\mathbb N}}
\newcommand{\zz}{{\mathbb Z}}

\newcommand\limproj{\mathop{\underleftarrow{\lim}}}
\newcommand\limind{\mathop{\underrightarrow{\lim}}}

\makeatother
\newtheorem{defin}{Definition}

\newtheorem{ex}{Example}
\newtheorem{rem}{Remark}
\theoremstyle{plain}
\newtheorem{prop}{Proposition}

\newtheorem{theo}{Theorem}
\newtheorem{lemma}{Lemma}
\newtheorem{corol}{Corollary}

\newcommand{\eqdef}{\stackrel{\rm def}{=}}

\newcommand{\Proof}{{\noindent \bf Proof}}
\newcommand{\lto}{\longrightarrow}
\newcommand{\Spec}{\mathop{ \rm Spec}}
\newcommand{\Spf}{\mathop{ \rm Spf}}

\newcommand{\Ker}{\mathop{\rm Ker}}
\newcommand{\Coker}{\mathop{\rm Coker}}

\newcommand{\Hom}{\mathop{\rm Hom}}
\newcommand{\End}{\mathop{\rm End}}
\renewcommand{\delta}{\partial}

\title{Formal groups arising from formal punctured ribbons.\footnote{the second and
the third authors are supported by RFBR grant no. 08-01-00095-a, by
INTAS grant 05-1000008-8118; besides the second author is supported
by grant of Leading Scientific Schools no. 1987.2008.1, by a program
of President of RF for supporting of young russian scientists (grant
no. MK-864.2008.1), and the third author is supported by grant of
Leading Scientific Schools no. 4578.2006.1, and by grant of National
Scientific Projects no. 2.1.1.7988.}}
\author{Herbert Kurke, \quad Denis Osipov, \quad Alexander Zheglov}

\date{}

\begin{document}

\maketitle

\begin{abstract} We investigate Picard functor of a formal punctured
ribbon. We prove that under some conditions this functor is
representable by a formal group scheme.

\vspace{0.5cm}


\vspace{0.5cm}


\end{abstract}

\section{Introduction.}

First we would like to give a motivation of our present research.

Let's consider a $\bc$-algebra $R$ with a derivation $\delta : R \to R$
$$
\delta(ab)= \delta(a)b + a \delta(b) , \quad a,b \in R \mbox{.}
$$
We construct a ring
$$R((\delta^{-1})) \quad : \quad \sum_{i \ll
+\infty} a_i \delta^{i}, a_i \in R$$
$$
[\delta, a] = \delta(a) , \quad \delta^{-1}a =a \delta^{-1} +
C_{-1}^{1} \delta(a) \delta^{-2} + C_{-1}^{2} \delta^2(a)
\delta^{-3} + \ldots \mbox{,}
$$
where $C_{i}^k$, $i \in \dz$, $k \in \dn$ is a binomial coefficient:
$$
C_i^k = \frac{i(i-1) \ldots (i-k+1)}{k(k-1)\ldots 1} \mbox{,} \quad
C_i^0=1 \mbox{.}
$$

Now we consider $R = \bc [[x]]$ with  usual derivation
$\delta(x)=1$. We add infinite number of "formal times" : $t_1, t_2,
\ldots$. There is a unique decomposition in the ring
$R((\delta^{-1}))[[t_1, t_2, \ldots]]$:
$$
\mbox{if} \qquad A \in R((\delta^{-1}))[[t_1, t_2, \ldots]] \mbox{,}
\qquad \mbox{then} \quad A=A_+ + A_-  \mbox{,}$$
$$\mbox{where}  \quad A_+ \in R[\delta][[t_1,
t_2, \ldots]], \quad A_- \in R[[\delta^{-1}]] \cdot \delta^{-1}
[[t_1, t_2, \ldots]] \mbox{.}
$$

Let $L \in R((\delta^{-1}))[[t_1, t_2, \ldots]]$ be of the following
type:
$$
L = \delta + a_1\delta^{-1} + a_2\delta^{-2} + \ldots , \qquad a_i
\in R[[t_1, t_2, \ldots]] \mbox{.}
$$

The classical {\bf KP}-hierarchy is the following infinite system of
equations, see~\cite{SW}:
$$
\frac{\partial L }{\partial t_n} = [(L^n)_+, L] \mbox{,} \qquad n
\in \nn \mbox{.}
$$
From this system it follows $$ \mbox{the {\bf KP} equation} \qquad
(4u_t-u'''-12uu')'=3u_{yy} \qquad \mbox{for}  \quad u(t,x,y)
\mbox{,}$$
$$ \mbox{and the {\bf KdV} equation} \qquad  4u_t-7u''' - 12uu'=0 \qquad \mbox{for} \quad
u(t,x)\mbox{.}$$ Solutions of {\bf KP}-hierarchy are obtained from
flows on Picard varieties of algebraic curves (for example,
solitons).

A.N.~Parshin gave in 1999 in~\cite{Pa0} the following generalization
of {\bf KP}-hierarchy. (A.~Zheglov modified it later in~\cite{Zhe}.)
Let $R = \bc[[x_1, x_2]] ((\delta_1^{-1}))$, where the derivation
$$\delta_1 \quad : \quad \bc[[x_1, x_2]] \to \bc[[x_1, x_2]], \qquad
\delta_1(x_1)=1, \quad \delta_1(x_2)=0 \mbox{.}$$ We consider a
derivation $\delta_2 : R \to R $, $\delta_2(x_1)=0$,
$\delta_2(x_2)=1$, $\delta_2(\delta_1)=0$. As before, we construct a
ring
$$
E=R((\delta_2^{-1})) = \bc[[x_1, x_2]]
((\delta_1^{-1}))((\delta_2^{-1})) \mbox{.}
$$
We add "formal times" $\{ t_k \}$, $k =(i,j) \in \zz \times \zz_+$. As
before, there is a decomposition (with respect to $\partial_2$):
$$E[[\{ t_k \}]]=
E_+[[\{ t_k \}]] \oplus E_-[[\{ t_k \}]] \mbox{.}$$ We consider $L,M
\in E[[\{ t_k \}]]$ such that
$$
L = \delta_1 + u_1 \delta_2^{-1} + \ldots  , \qquad M = \delta_2 +
v_1\delta_2^{-1} + \ldots\mbox{,} $$
$$\mbox{where} \quad u_i, v_i
\in R[[\{ t_k \}]] \mbox{.}
$$
Let $N = (L,M)$ and $[L,M]=0$, then {\it hierarchy} is
$$
\frac{\partial N}{\partial t_k} = V_N^k \mbox{,}
$$
$$
\mbox{where} \quad V_N^k= ([(L^i M^j)_+, L], [(L^iM^j)_+, M])
\mbox{,}
$$
$$
 k=(i,j) \in \zz \times \zz_+ , \quad i \le \alpha j, \quad \alpha > 0 \mbox{.}
$$

There is the following property. Let $L,M \in E$ such that they
satisfy conditions for Parshin's hierarchy when all the times $t_k
=0$. Then there is $S \in 1 + E_-  \subset E$ such that
$$
L= S^{-1} \delta_1 S    ,  \qquad M= S^{-1} \delta_2 S \mbox{.}
$$

Besides, the ring $E$ acts $\bc$-linearly on $\bc((u))((t))$ (and on
the set of $\bc$-vector subspaces of $\bc((u))((t))$) in the
following way:
$$
E/ E \cdot (x_1,x_2) = \bc((u))((t)), \qquad \delta_1^{-1} \mapsto
u, \quad \delta_2^{-1} \mapsto t,
$$
now $E$ acts naturally on the left on $E/ E \cdot (x_1,x_2)$.

In classical {\bf KP}-hierarchy an analogous action is an action of
the ring of pseudodifferential operators on the set of Fredholm
subspaces of $\bc((t))$, or more generally, on the Sato Grassmanian.
This action gives the flows on generalized Jacobians of algebraic
curves, which are solutions of~{\bf KP}-hierarchy,
see~\cite[\S1]{Pa}.

In article~\cite{Ku} we investigated new geometric objects
$\xo_{\infty}=(C, \ca)$, which are ringed spaces: formal punctured
ribbons with the underlying topological space $C$ as an algebraic
curve. (For simplicity we call such objects "ribbons".) Examples of
ribbons come from Cartier divisors on algebraic surfaces.

We are working in formal algebraic language, therefore originally we
assume that ribbons are defined over any ground field $k$. But in
many places of this article we will additionally assume that $k$ is
an algebraically closed field of characteristic zero.

We introduced the notion of a torsion free sheaf on a ribbon. An
importance of such sheaves followed from theorem~1 of
article~\cite{Ku}, where we proved that torsion free sheaves on some
ribbons plus some geometrical data such as formal trivialization of
sheaves, local parameters at smooth points of ribbons and so on are
in one-to-one correspondence with {\it generalized Fredholm
subspaces} of two-dimensional local field $k((u))((t))$ (see also
section~\ref{tfs} of this article).

In this article we investigate torsion free sheaves on ribbons
$(C,\ca)$ and proved  that if the underlying curve $C$ of a ribbon
is a smooth curve and for any small open $U \subset C$ there are
sections $a \in \Gamma(U, \ca_1)$, $a^{-1} \in \Gamma(U, \ca_{-1})$,
then every torsion free sheaf on the ribbon $(C, \ca)$ is a locally
free sheaf on a ringed space $(C, \ca)$, see
proposition~\ref{freeness}. We remark that this condition is
satisfied, for example, when the ribbon $(C, \ca)$ comes from a
smooth curve $C$, which is a Cartier divisor on algebraic surface.

Therefore it is important to study locally free sheaves on ribbons
$(C, \ca)$. We restrict ourself to the Picard group of a ribbon.
In~\cite{Ku} we investigated the Picard group as a set, see
proposition~5 and example~8 in~\cite{Ku}. But it was not clear, what
are the deformations (local or global) of elements of
$Pic(\xo_{\infty})$. We study  the groups $Pic(\xo_{\infty,S})$ and
$Pic(X_{\infty, S})$ for an arbitrary affine scheme $S$ as functors
$Pic_{\sxo_{\infty}}$ and $Pic_{X_{\infty}}$ on the category of
affine schemes from the point of view of representability or formal
representability of these functors by a scheme or a formal scheme,
see section~\ref{fpfb} and section~\ref{picfunrib}. We note that
the  functor $Pic_{X_{\infty}}$ is mapped in the
 functor $Pic_{\sxo_{\infty}}$.

At first, we study  the tangent spaces to these Picard functors. In
article~\cite{ZhO} the "picture cohomology" $\h^0(W)$, $\h^1(W)$,
$\h^2(W)$ were introduced for a generalized Fredholm subspace $W$ of
a two-dimensional local field. These cohomology groups coincide with
the cohomology groups of a line bundle on an algebraic surface when
a ribbon and a line bundle on it come from an algebraic surface and
a line bundle on this surface. In section~\ref{piccoh} we
investigate the picture cohomology groups of generalized Fredholm
subspaces $W$ and related them with some groups which depend on
cohomology groups of sheaves $\cf_{W}$ and $\cf_{W, 0}$ on the curve
$C$, where $W \longleftrightarrow (\cf_W, \ldots )$ is a generalized
Krichever-Parshin correspondence from~\cite[\S5]{Ku}. Due to this
result we obtained that the kernel of the natural map from tangent
space of functor $Pic_{X_{\infty}}$ to tangent space of functor
$Pic_{\sxo_{\infty}}$  coincides with the first picture cohomology
of the structure sheaf, and cokernel of this map coincides with the
second picture cohomology of the structure sheaf, see
proposition~\ref{tspc}.

Further, in section~\ref{fpfb}, we investigate the Picard functors
on ringed spaces $X_{\infty}$ and $\xo_{\infty}$ as formal functors
on Artinian rings. We prove that if the first picture cohomology
group of the structure sheaf of a ribbon $\xo_{\infty}$ is
finite-dimensional over the ground field $k$ and $\chara  k = 0$,
then the formal Picard functor $\widehat{Pic}_{\sxo_{\infty}}$ is
representable by a formal group, which can be decomposed in the
product of two formal groups, where the first one is connected with
the formal Picard functor $\widehat{Pic}_{X_{\infty}}$ and the
second one coincides with the formal Brauer group of algebraic
surface when the ribbon $\xo_{\infty}$  comes from an algebraic
surface and a curve on it (see corollary~\ref{cor1} of
proposition~\ref{prop4}).

In section~\ref{picfunrib} we prove (under some condition) the
global representability of the Zariski sheaf associated  with the
presheaf (or functor) $\widetilde{Pic}_{\sxo_{\infty}}$ on the
category of Noetherian schemes. The condition we mean here is equivalent
to $\h^1(W) =0$, where $W \subset k((u))((t))$ correspond to the
structure sheaf of the ribbon $\xo_{\infty}$ plus some local
parameters via the Krichever-Parshin map. The functor
$\widetilde{Pic}_{\sxo_{\infty}}$ classifies invertible sheaves plus
morphisms of order on the ribbon $\xo_{\infty}$, see
section~\ref{rpf}. Then we prove in theorem~\ref{represent} that
(noncanonically) Picard scheme of $\xo_{\infty}$ is the product of the
Picard scheme of $X_{\infty}$ (see section~\ref{repxinf}) and the formal Brauer
group of $\xo_{\infty}$ (see section~\ref{bgr}) when the ground
field $k$ is an algebraically closed field and $\chara k =0$.
From this result and under the same conditions we obtain in theorem~\ref{theolast}
the
global representability of the Zariski sheaf associated  with the
presheaf  ${Pic}_{\sxo_{\infty}}$ on the
category of Noetherian schemes.

At last, let's note here that there are many activities in the
direction of constructing of geometric objects which encode spectral
properties of commutative rings of germs of differential operators
in the 2-dimensional case. For an (incomplete) survey on recent
activities one can consult~\cite{P}.

This work was partially done during our stay at the Erwin
Schr\"{o}dinger International Institute for Mathematical Physics
from January 8 - February 2, 2008 and at the Mathematisches
Forschungsinstitut Oberwolfach during a stay within the Research in
Pairs Programme from January 25 - February 7, 2009. We want to thank
these institutes for hospitality and excellent working conditions.
We are grateful also to A.N. Parshin for some discussions.

\section{Torsion free sheaves on ribbons.}
\label{tfs}

We recall the general definition of ribbon from~\cite{Ku}.

Let $S$ be a Noetherian base scheme.
\begin{defin}[\cite{Ku}] \label{ribbon}
A ribbon $(C, \ca)$ over $S$ is given by the following data.
\begin{enumerate}
\item\label{a}
 A flat family of reduced algebraic curves $\tau :  C\rightarrow S$.
\item \label{b} A sheaf $\ca$ of commutative $\tau^{-1}\co_S$-algebras
on $C$.
\item \label{c} A descending sheaf filtration $(\ca_i)_{i\in \sdz}$ of $\ca$ by
$\tau^{-1}\co_S$-submodules which satisfies  the following axioms:
\begin{enumerate}
\item  \label{i} $\ca_i\ca_j \subset \ca_{i+j}$, $1\in \ca_0$ (thus $\ca_0$
is a subring, and for any $i \in \dz$ the sheaf $\ca_i$ is a
$\ca_0$-submodule);
\item \label{ii}
$\ca_0/\ca_1$ is the structure sheaf $\co_C$ of $C$;
\item
\label{iii} for each $i$ the sheaf $\ca_i/\ca_{i+1}$ (which is a
$\ca_0/\ca_1$-module by (\ref{i})) is a coherent sheaf on $C$, flat
over $S$, and for any $s \in S$ the sheaf $\ca_i/\ca_{i+1}
\mid_{C_s}$ has no coherent subsheaf with finite support, and is
isomorphic to $\co_{C_S}$ on a dense open set;
\item
\label{iv}  $\ca =\limind\limits_{i \in \sdz} \ca_i$, and
$\ca_i=\limproj\limits_{j>0}\ca_i/\ca_{i+j}$ for each $i$.
\end{enumerate}
\end{enumerate}
\end{defin}

Sometimes we shall denote a ribbon $(C, \ca)$ over $\Spec R$, where
$R$ is a ring, as a ribbon over $R$.

There is the following example of a ribbon.  If $X$ is an algebraic
surface over a field $k$, and $C\subset X$ is a reduced effective
Cartier divisor, we obtain a ribbon $(C, \ca )$ over $k$, where
$$ \ca := \co_{\hat{X}_C}(*C)=\limind\limits_{i \in \sdz} \co_{\hat{X}_C}(-iC)
= \limind\limits_{i \in \sdz} \: \limproj\limits_{j \ge 0} J^{i}/
J^{i+j}
$$
$$
 \ca_i := \co_{\hat{X}_C}(-iC) = \limproj\limits_{j \ge 0} J^{i}/
J^{i+j}  \mbox{,} \quad i \in \dz  \mbox{,}$$
 where $\hat{X}_C$ is the formal scheme which is the completion of $X$ at
 $C$,
 and $J$ is the ideal sheaf of $C$ on $X$ (the sheaf $J$ is an invertible
 sheaf). We shall say that a ribbon which is constructed by this example
 is {\em "a ribbon which comes from an algebraic surface"}.

We recall the definition of a torsion free sheaf on a ribbon
from~\cite{Ku}.
\begin{defin}[\cite{Ku}] \label{tfsh}
Let $\xo_{\infty}=(C,\ca )$ be a ribbon over a scheme $S$. We say
that $\cn$ is a torsion free sheaf of rank $r$ on $\xo_{\infty}$ if
  $\cn$ is a sheaf of $\ca$-modules on $C$ with a descending
filtration $(\cn_i)_{i\in \sdz}$ of $\cn$ by $\ca_0$-submodules
which satisfies the following axioms.
\begin{enumerate}
\item
 $\cn_i\ca_j\subseteq \cn_{i+j}$ for any $i,j$.
\item \label{ittf}
For each $i$ the sheaf $\cn_i/\cn_{i+1}$ is a coherent sheaf on $C$,
flat over $S$, and for any $s \in S$  the sheaf
$\cn_i/\cn_{i+1}|_{C_S}$ has no coherent subsheaf with finite
support, and is isomorphic to $\co_{C_S}^{\oplus r}$ on a dense open
set.
\item
 $\cn =\limind\limits_i \cn_i$ and
$\cn_i=\limproj\limits_{j>0}\cn_i/\cn_{i+j}$ for each $i$.
\end{enumerate}
\end{defin}

\begin{rem}
\label{flat}
Note that the sheaf $\cn$ is flat over $S$. To show this note that all sheaves $\cn_i/\cn_{i+j}$ are, clearly, flat over $S$ (see, e.g. \cite{Ku}, prop.1). So, for any ideal sheaf $J$ on $S$ we have the embeddings $0\rightarrow \tau^*J\otimes \cn_i/\cn_{i+j} \rightarrow \cn_i/\cn_{i+j}$ by the flatness criterium for modules. This imply that we have embeddings  $0\rightarrow \tau^*J\otimes \cn_i \rightarrow \cn_i$ for any $i$ and embeddings $0\rightarrow \tau^*J\otimes \cn \rightarrow \cn $. Therefore, $\cn$ is flat over $S$.
\end{rem}

There is the following example of a torsion free sheaf of rank $r$
on a ribbon $\xo_{\infty}=(C, \ca)$ which comes from an algebraic
surface $X$. Let $E$ be a locally free sheaf of rank $r$
 on the surface $X$. Then
$$\eo_C:=\limind\limits_i\limproj\limits_jE(iC)/E(jC)$$
is a torsion free sheaf of rank $r$ on $\xo_{\infty}$. We shall say
that a torsion free sheaf on a ribbon constructed after this example
 is {\em "a sheaf which comes from a locally free sheaf on an algebraic surface"}.

In~\cite{Ku} we defined the notion of a smooth point of a ribbon,
the notion of formal local parameters at a smooth point of a ribbon,
and the notion of a smooth point of a torsion free sheaf on a
ribbon, see definitions 9, 10, 12 from~\cite{Ku}. We remark that
these notions coincide with the usual notions (i.e. used in
\cite{Pa}, \cite{Os}) when a ribbon comes from algebraic surface, a
torsion free sheaf on a ribbon comes from a locally free sheaf on
this surface and so on.

Let $k$ be a field.  We recall (see, for example,~\cite{ZhO}) that a
$k$-subspace $W$ in $k((u))^{\oplus r}$ is called a {\em Fredholm}
subspace if
$$
\dim_k W \cap k[[u]]^{\oplus r} < \infty  \qquad \mbox{and }\qquad
\dim_k \frac{k((u))^{\oplus r}}{W + k[[u]]^{\oplus r}} < \infty
\mbox{.}
$$

\begin{defin}[cf.~definitions 15--16 of~\cite{Ku}]
For a $k$-subspace $W$ in $k((u))((t))^{\oplus r}$, for  $n \in \dz$
let
$$
W(n) = \frac{W \cap t^nk((u))[[t]]^{\oplus r}}{W \cap
t^{n+1}k((u))[[t]]^{\oplus r}}
$$
be a $k$-subspace in $k((u))^{\oplus r} =
\frac{t^nk((u))[[t]]^{\oplus r}}{t^{n+1}k((u))[[t]]^{\oplus r}}$.

A closed $k$-subspace $W$ in $k((u))((t))^{\oplus r}$ is called {\em a
generalized Fredholm subspace} iff for any $n \in \dz$ the
$k$-subspace $W(n)$ in $k((u))^{\oplus r}$ is a Fredholm subspace.
\end{defin}

The following definition is from~\cite{Ku}.
\begin{defin}[\cite{Ku}]
Let a $k$-subalgebra $A$ in $k((u))((t))$ be a generalized Fredholm
subspace. Let a $k$-subspace $W$ in $k((u))((t))^{\oplus r}$ be a
generalized Fredholm subspace. We say that  $(A,W)$ is {\em a Schur
pair} if $A \cdot W \subset W$.
\end{defin}

Now we recall the main theorem of~\cite{Ku}.
\begin{theo}[\cite{Ku}]
\label{th1}  Schur pairs $(A,W)$ from $k((u))((t)) \oplus
k((u))((t))^{\oplus r}$ are in one-to-one correspondence with the data
$(C,\ca , \cn ,P, u,t, e_P)$  up to an isomorphism, where  $C$ is a
projective irreducible curve over a field $k$, $(C, \ca)$ is a
ribbon, $\cn$ is a torsion free sheaf of rank $r$ on this ribbon,
$P$ is a point of $C$ which is a smooth point of $\cn$, $u,t$ are
formal local parameters of this ribbon at $P$, $e_P$ is a formal
local trivialization of $\cn$ at $P$.
\end{theo}

The goal of this section is to show that in "good" cases a torsion
free sheaf  on a ribbon $(C, \ca)$ is indeed a locally free sheaf on
this ribbon, that is an element of the set $\check{H}^1(C,
GL_r(\ca))$.

We recall the following condition from~\cite{Ku} (see lemma~4
of~\cite{Ku}).
\begin{defin}
\label{uslovie} The sheaf $\ca$ of a ribbon $(C, \ca )$ satisfies
(**) if the following condition holds: there is an affine open
cover $\{U_{\alpha}\}_{\alpha\in I}$ of $C$ such that for any
$\alpha \in I$  there is  an invertible section $a\in \ca_1
(U_{\alpha})  \subset \ca (U_{\alpha})$ such that $a^{-1} \in
\ca_{-1}(U_{\alpha})$.
\end{defin}
We note that when a ribbon $(C, \ca)$ comes from an algebraic
surface, then the condition (**) is satisfied for $\ca$. In this
case elements $a$ (from definition~\ref{uslovie}) come from local
equations of the Cartier divisor $C$ on the algebraic surface.

We have the following proposition.
\begin{prop}
\label{freeness}
Let the sheaf $\ca$ of a ribbon $(C, \ca)$ over $S$ satisfy
the condition $(**)$. Let a torsion free sheaf $\cn$ of rank $r$ on
the ribbon $(C,\ca)$ satisfy the following condition: the sheaf
$\cn_0/\cn_1$ is a locally free sheaf on $C$. Then the
sheaf $\cn$ is a locally free sheaf of rank $r$ on the ribbon $(C,
\ca)$.
\end{prop}
\begin{rem}
\label{rem2} The condition on the sheaf $\cn_0/\cn_1$ from
proposition is satisfied, for example, if $(C, \ca )$ is a ribbon
over a field $k$ with $C$ a smooth curve, since from
definition~\ref{tfsh} we have that $\cn_0/\cn_1$ is a torsion free
sheaf. It is also satisfied if $C$ is a flat family of smooth curves
over $S$ due to the following fact: if $\cf$ is a coherent sheaf on
$C$ which is flat over $S$, and such that $\cf |_{C\times_S s}$ is a
locally free sheaf for any $s\in S$, then $\cf$ is a locally free
sheaf on $C$. (The last statement follows easily from  Nakayama's
lemma)
\end{rem}
\Proof. We shall prove that if an open affine $V$ of $C$ such that
$V \subset U_{\alpha}$ for some $\alpha \in I$ (see
definition~\ref{uslovie}) and $\cn_0/ \cn_1 \mid_V \simeq
\co_V^{\oplus r}$, then $\cn \mid_V \simeq \ca^{\oplus r} \mid_V$.

Let $\bar{c}_1, \ldots, \bar{c}_r \in \cn_0 / \cn_1 (V)$ be a basis
over $\co_C (V)$. We choose some elements $c_1, \ldots, c_r \in
\cn_0(V)$ such that for any $1 \le  i \le r$ the element $c_i$ maps
to the element $\bar{c}_i$ under the natural map $\cn_0(V) \lto
\cn_0/\cn_1(V)$. (By proposition~3 of~\cite{Ku}, the last map is a
surjective map, therefore such elements $c_1, \ldots, c_r$ exist.)

We consider a map of $\ca \mid_V$-modules:
$$
\phi : \ca^{\oplus r} \mid_V  \lto \cn \mid_V    \qquad \qquad
\phi(\bigoplus\limits_{1 \le i \le r} a_i)= \sum\limits_{1 \le i \le
r} a_i \cdot c_i \mbox{,}
$$
where $a_i \in \ca(U)$ for $1 \le i \le r$, an open $U \subset V$.

At first, we show that the map $\phi$ is a surjective map of
sheaves. Let an element $b \in \cn(U)$, $b \ne 0$ for an open $U
\subset V$. Then $b \in \cn_{l_1}(U) \setminus \cn_{l_1+1}(U)$ for
some $l_1 \in \dz$. Therefore an element $b_1=a^{-l_1} \cdot b \in
\cn_0(U) \setminus \cn_1(U)$, where $a \in \ca_1(V) \setminus
\ca_2(V)$ such that $a^{-1} \in \ca^{-1}(V)$. Let $\bar{b}_1 \in
\cn_0/\cn_1(U)$ be the image of the element $b_1$. We have
$$
\bar{b}_1 = \sum_{1 \le i \le r} \bar{e}_{1,i} \cdot \bar{c}_i
\mbox{,}
$$
where $\bar{e}_{1,i} \in \co_C(U)$, $1 \le i \le r$. We choose some
elements $e_{1,i} \in \ca_0(U)$ such that for any $1 \le i \le r$
the image of the element $e_{1,i}$ in $\ca_0(U) / \ca_1(U) =
\co_C(U)$ coincides with the element $\bar{e}_{1,i}$ (see also
proposition~3 of~\cite{Ku}). Now if $b_1 \ne \sum_{1 \le i \le r}
e_{1,i} \cdot c_i$ , then an element
$$
(b_1 - \sum_{1 \le i \le r} e_{1,i} \cdot c_i) \quad \in \quad
\cn_{l_2}(U) \setminus \cn_{l_2 +1}(U)
$$
for some $l_2 \in \dn$, where $l_2 \ge 1$. Therefore an element
$$
b_2= a^{-l_2} \cdot (b_1 - \sum_{1 \le i \le r} e_{1,i} \cdot c_i)
\quad \in \quad \cn_0(U) \setminus \cn_1(U) \mbox{.}
$$
And we can repeat the same procedure with $b_2$ as with $b_1$
before, and so on.

Now an element
$$
d = a^{l_1} \cdot ( \bigoplus_{1 \le i \le r}  e_{1,i}  + a^{l_2}
\cdot ( \bigoplus_{1 \le i \le r} e_{2,i}   + \ldots) )
$$
is well defined in $\ca(U)^{\oplus r}$ as a convergent infinite
series, because $\ca(U)$ is a complete space and $l_n \ge 1$ for
$n>1$. And, by construction, $\phi(d) = b$, because $\ca(V)^{\oplus
r}$ is a Hausdorff space. Therefore $\phi$ is a surjective map.

Second, we show that $\phi$ is an injective map of $\ca^{\oplus r}
\mid_V$-modules. Let the sheaf $\ck$ be a kernel of the map $\phi$.
Let $g \in \ck(U)$, $g \ne 0$ for some open $U \subset V$. We have
$g \in \ca_l(U)^{\oplus r} \setminus \ca_{l+1}(U)^{\oplus r} $ for
some $l \in \dz$, then $a^{-l} \cdot g \in \ca_0(U)^{\oplus r}
\setminus \ca_{1}(U)^{\oplus r}$. Let $e=(e_1, \ldots, e_r) \in
\co_C(U)^{\oplus r}$ be the image of $a^{-l} \cdot g$ under the
natural map. Since $a^{-l} \cdot g \in \ck$, we have
$$
\sum_{1 \le i \le r} e_i \cdot \bar{c}_i = 0 \mbox{.}
$$
Therefore for any $1 \le i \le r$ $e_i=0$, because $\bar{c}_1,
\ldots, \bar{c}_r$ is a basis over $\co_C(U)$. Hence, $a^{-l} \cdot
g \in \ca_1(U)^{\oplus r}$. We have a contradiction.
\begin{flushright}
$\square$
\end{flushright}

\section{"Picture cohomology".}
\label{piccoh}
Let $k$ be a field. Let  $W$ be a  $k$-subspace in
$k((u))((t))^{\oplus r}$. Let
$$
\co_1 = k((u))[[t]] \qquad \mbox{,} \qquad \co_2=k[[u]]((t))
 $$
be  $k$-subspaces in $k((u))((t))$. We consider the following
complex.
\begin{equation} \label{rescom}
(W \cap \co_2^{\oplus r}) \oplus (W \cap \co_1^{\oplus r}) \oplus
(\co_1^{\oplus r} \cap \co_2^{\oplus r})  \lto
 W \oplus  \co_2^{\oplus r} \oplus
   \co_1^{\oplus r}
\lto    k((u))((t))^{\oplus r}
\end{equation}
where the first map is given by
$$
(a_0, a_1, a_2)            \mapsto (a_1 - a_0, a_2 - a_0, a_2 -a_1)
$$
and the second by
$$
(a_{01}, a_{02}, a_{12})
               \mapsto  a_{01} - a_{02} + a_{12} \mbox{.}
$$

\begin{rem} \label{r2}
We suppose that a $k$-subspace $W \subset k((u))((t))^{\oplus r}$ is
a part of a Schur pair
$$(A,W)
\subset k((u))((t)) \oplus k((u))((t))^{\oplus r} \mbox{.}$$
 Let, by
theorem~\ref{th1}, the pair $(A,W)$ correspond to the data $(C,\ca ,
\cn ,P, u,t, e_P)$. We suppose that the ribbon $(C, \ca)$ comes from
an algebraic projective surface $X$, the torsion free sheaf $\cn$
comes from a locally free sheaf $\cf$ on $X$ and so on. It means
that the data $(C,\ca , \cn ,P, u,t, e_P)$ comes from the data $(X,
C, \cf, P, u, t, e_P)$, where $X$ is an algebraic projective
surface, $C$ is a reduced effective Cartier divisor, $\cf$ is a
locally free sheaf of rang $r$ on $X$, $P \in C$ is a point which is
a smooth point on $X$ and $C$, $u,t$ are formal local parameters of
$X$ at $P$ such that $t=0$ gives the curve $C$ on $X$ in a formal
neigbourhood of $P$ on $X$, $e_P$ is a formal trivialization of
$\cf$ at $P$. We suppose also that $X$ is a Cohen-Macaulay surface
and $C$ is an ample divisor on $X$. Then it was proved in~\cite{Os,
Pa} that the cohomology groups of
 complex~(\ref{rescom}) coincide with the cohomology groups $H^*(X, \cf)$.
\end{rem}

The goal of this section is to relate in general situation the
cohomology groups of complex~(\ref{rescom}) with the cohomology
groups of sheaves $\cn_i$, where $\cn$ is a torsion free sheaf on
the ribbon $(C, \ca)$ when, for example, this ribbon does not come
from an algebraic surface.

\begin{lemma} \label{l1}
Let $W$ be a $k$-subspace in $k((u))((t))^{\oplus r}$. Then the
cohomology groups of complex~(\ref{rescom}) coincide with the
following $k$-vector spaces:
$$
\h^0(W)= W \cap \co_1^{\oplus r} \cap \co_2^{\oplus r} \mbox{,}
$$
$$
\h^1(W)=\frac{W \cap (\co_1^{\oplus r} + \co_2^{\oplus r})} { W \cap
\co_1^{\oplus r} + W \cap \co_2^{\oplus r}} \mbox{,}
$$
$$
\h^2(W)= \frac{k((u))((t))^{\oplus r}}{ W + \co_1^{\oplus r} +
\co_2^{\oplus r}} \mbox{.}
$$
\end{lemma}
\Proof. We have the following exact sequence:
$$
0 \lto  \co_1^{\oplus r} \cap \co_2^{\oplus r}  \lto \co_1^{\oplus
r} \oplus \co_2^{\oplus r} \lto \co_1^{\oplus r} + \co_2^{\oplus r}
\lto 0 \mbox{,}
$$
where $\co_1^{\oplus r} + \co_2^{\oplus r}$ is considered as a
$k$-subspace in $k((u))((t))^{\oplus r}$. Now we take the
factor-complex of complex~(\ref{rescom}) by the following acyclic
complex:
$$
\co_1^{\oplus r} \cap \co_2^{\oplus r} \lto \co_1^{\oplus r} \cap
\co_2^{\oplus r} \lto 0 \mbox{.}
$$

We obtain the following complex:
$$
(W \cap \co_2^{\oplus r}) \oplus (W \cap \co_1^{\oplus r}) \lto
 W \oplus  (\co_2^{\oplus r} +
   \co_1^{\oplus r})
\lto    k((u))((t))^{\oplus r} \mbox{.}
$$
The cohomology groups of the last complex coincide with the
cohomology groups of complex~(\ref{rescom}). Therefore the statement
of this lemma is evident now.
\begin{flushright}
$\square$
\end{flushright}

\begin{defin}
The $k$-vector spaces $\h^i(W)$, $0 \le i \le 2$ are called {\em
"the picture cohomology"} of $W$.
\end{defin}

We have the following theorem.
\begin{theo} \label{pc}
Let $W$ be a $k$-subspace in $k((u))((t))^{\oplus r}$ such that the
$k$-space $W$ is a part of a Schur pair
$$(A,W)
\subset k((u))((t)) \oplus k((u))((t))^{\oplus r} \mbox{.}$$
 Let, by
theorem~\ref{th1}, the pair $(A,W)$ correspond to the data $(C,\ca ,
\cn ,P, u,t, e_P)$ (see the formulation of theorem~\ref{th1}). Then
\begin{equation} \label{eq1}
\h^0(W) = H^0(C, \cn_0)  \mbox{,}
\end{equation}
\begin{equation} \label{eq2}
\h^1(W) = \frac{H^0(C, \cn/\cn_0)}{\frac{H^0(C, \cn)}{H^0(C,
\cn_0)}} \mbox{,}
\end{equation}
\begin{equation}  \label{eq3}
\h^2(W) = H^1(C, \cn/\cn_0) \mbox{.}
\end{equation}
\end{theo}
\begin{rem}
Here and further in the article we consider cohomology in Zariski
topology, if another topology is not specified.
\end{rem}
\Proof. By definition of a torsion free sheaf on a ribbon, we have
$\cn_0 = \limproj\limits_{i > 0} \cn_0/\cn_i$. Therefore
$$
H^0(C, \cn_0)=\limproj\limits_{i > 0} H^0(C,\cn_0/\cn_i)=
\limproj\limits_{i > 0} (W(0,i) \cap \co_2^{\oplus r} ) = W \cap
\co_1^{\oplus r} \cap \co_2^{\oplus r} \mbox{,}
$$
where $W(l,i) \eqdef \frac{W \cap t^l \cdot \co_1^{\oplus r}}{W \cap
t^i \cdot
 \co_1^{\oplus r}}$, $l < i \in \dz$. Here we used theorem~2 of~\cite{Os}, where the complex
 was constructed  which calculates in our case the cohomology
 groups of the coherent sheaf $\cn_0/\cn_i$ of $\ca_0/\ca_i$-modules
 on the $1$-dimensional scheme $X_{i-1}=(C, \ca_0/\ca_i)$. Formula~(\ref{eq1}) is proved.

Now we will prove formula~(\ref{eq2}). We have
$$
\h^1(W)= \frac{W \cap (\co_1^{\oplus r} + \co_2^{\oplus r}) }{W \cap
\co_1^{\oplus r} + W \cap \co_2^{\oplus r}} = \frac{\frac{W \cap
(\co_1^{\oplus r} + \co_2^{\oplus r})}{W \cap \co_1^{\oplus
r}}}{\frac{W \cap \co_1^{\oplus r} + W \cap \co_2^{\oplus r}}{W \cap
\co_1^{\oplus r}}} = \frac{\frac{W \cap (\co_1^{\oplus r} +
\co_2^{\oplus r})}{W \cap \co_1^{\oplus r}}}{\frac{W \cap
\co_2^{\oplus r}}{W \cap \co_2^{\oplus r} \cap \co_1^{\oplus r}}}
\mbox{.}
$$
We note that we have
$$
W \cap \co_2^{\oplus r} = \limind\limits_i \limproj\limits_{j > i}
W(i,j) \cap \co_2^{\oplus r} = \limind\limits_i \limproj\limits_{j >
i} H^0(C, \cn_i / \cn_j) = H^0(C, \cn) \mbox{.}
$$
Here we used theorem~2 of~\cite{Os} for the coherent sheaf
$\cn_i/\cn_j$ of $\ca_0/\ca_{j-i}$-modules
 on the $1$-dimensional scheme $X_{j-i-1}=(C, \ca_0/\ca_{j-i})$.
Therefore, using it and formula~(\ref{eq1}), we have
$$
\frac{W \cap \co_2^{\oplus r}}{W \cap  \co_2^{\oplus r} \cap
\co_1^{\oplus r}} = \frac{H^0(C, \cn)}{H^0(C, \cn_0)} \mbox{.}
$$
Hence, to prove formula~(\ref{eq2}), we have to check that
\begin{equation} \label{eq4}
\frac{W \cap (\co_1^{\oplus r} + \co_2^{\oplus r})}{W \cap
\co_1^{\oplus r}} = H^0(C, \cn/ \cn_0) \mbox{.}
\end{equation}

By proposition~3 of~\cite{Ku} we have that $H^1(C \setminus p,
\cn_0) =0$. Therefore from the exact triple of sheaves on the curve
$C$
$$
0 \lto \cn_0 \lto \cn \lto \cn/\cn_0 \lto 0
$$
we have
$$
H^0(C \setminus p, \cn / \cn_0) = \frac{H^0(C \setminus p,
\cn)}{H^0(C \setminus p, \cn_0)}= \frac{W}{W \cap \co_1^{\oplus r}}=
\frac{W + \co_1^{\oplus r}}{\co_1^{\oplus r}} \mbox{.}
$$
Now, as an inductive limit of complexes from theorem~2 of~\cite{Os},
we obtain that
$$
H^0(C, \cn/ \cn_0) = \frac{W + \co_1^{\oplus r}}{\co_1^{\oplus r}}
\cap \frac{\co_1^{\oplus r} + \co_2^{\oplus r}}{\co_1^{\oplus r}}
\mbox{,}
$$
where the intersection is considered in the $k$-vector space
$\frac{k((u))((t))^{\oplus r}}{\co_1^{\oplus r}}$.

There is a natural isomorphism of $k$-subspaces in the $k$-vector
space $\frac{k((u))((t))^{\oplus r}}{\co_1^{\oplus r}}$:
$$
\frac{W \cap (\co_1^{\oplus r} + \co_2^{\oplus r})}{W \cap
\co_1^{\oplus r}} = \frac{W }{W \cap \co_1^{\oplus r}} \cap
\frac{\co_1^{\oplus r} + \co_2^{\oplus r}}{\co_1^{\oplus r}}
\mbox{.}
$$
Therefore we checked formula~(\ref{eq4}). Hence, we proved
formula~(\ref{eq2}).

Now we will prove  formula~(\ref{eq3}). We have
$$
H^1(C, \cn/ \cn_0) = \limind\limits_{i < 0 } H^1(C, \cn_i / \cn_0)=
\limind\limits_{i < 0 } \frac{\frac{t^{i} \cdot \co_1^{\oplus
r}}{\co_1^{\oplus r}}}{ \frac{W \cap t^i \cdot \co_1^{\oplus
r}}{\co_1^{\oplus r}} + \frac{\co_2^{\oplus r}}{\co_2^{\oplus r}
\cap \co_1^{\oplus r}} } = \frac{k((u))((t))^{\oplus r}}{W +
\co_1^{\oplus r} + \co_2^{\oplus r}} \mbox{.}
$$
Here we used the complex from theorem~2 of~\cite{Os} to calculate
the first cohomology group of coherent sheaf $\cn_i / \cn_0$ on the
scheme $X_{-i-1}= (C, \ca_0 / \ca_{-i})$, $i < 0$.
Formula~(\ref{eq3}) is proved.
\begin{flushright}
$\square$
\end{flushright}

\section{Formal Picard group and formal Brauer group.}
\label{fpfb}
\subsection{Picard functor}

Let $Y \lto Z$ and $X \lto Z$ be morphisms of schemes. We consider
 $W= Y \times_Z X$ with the natural projection maps $p : W
\lto Y$ and $q : W \lto Z$. Let $\cf$ be an $\co_Y$-module sheaf on
$Y$, and $\cg$ be an $\co_X$-module sheaf on $X$. We recall the
definition of $\co_{W}$-module sheaf $\cf \boxtimes_{\co_Z} \cg $ on
$W$:
$$
\cf \boxtimes_{\co_Z} \cg  \eqdef p^*(\cf) \otimes_{\co_W} q^*(\cg)
\mbox{.}
$$

Now we recall the definition of a base change for a ribbon from
section~2.2 of~\cite{Ku}.

For a ribbon $\xo_{\infty} = (C, \ca)$ over $S$, and a morphism
$\alpha: S' \longrightarrow S$ of Noetherian schemes we define a
base change ribbon $\xo_{\infty, S'} = (C_{S'}, \ca_{S'})$ over $S'$
in the following way:
$$
C_{S'} :=C \times_S S'  \mbox{,}
$$
$$
\ca_{S', j} :=\limproj\limits_{i\ge 1}
(\ca_j/\ca_{j+i})\boxtimes_{\co_S}\co_{S'}, \mbox{\quad} \ca_{S'}:=\limind_j\ca_{S', j}
$$
for any $j \in \dz$. Sometimes we will denote $\ca_{S', j}$ by
$\ca_j\widehat{\boxtimes}_{\co_S} \co_{S'}$, and $\ca_{S'}$ by
$\ca\widehat{\boxtimes}_{\co_S} \co_{S'}$.

By the ribbon $\xo_{\infty} = (C, \ca)$ over $S$ we construct a
locally ringed space $ X_{\infty}= (C, \ca_0)$ over $S$. And also we
define a base change locally ringed space $ X_{\infty, S'}= (C_{S'},
\ca_{S', 0})$ for the morphism $\alpha: S' \longrightarrow S$ of
Noetherian schemes.

\begin{rem}
\label{locfin}
Let $\cm$ be an $\co_{S'}$-module sheaf on $S'$. Then we construct
the sheaf of $\ca_{S',j}$-modules $\ca_j\widehat{\boxtimes}_{\co_S}
\cm := \limproj\limits_{i\ge 1} (\ca_j/\ca_{j+i}) \boxtimes_{\co_S}
\cm$ on $C_{S'}$ for any $j \in \dz$, and the sheaf of
$\ca_{S'}$-modules $\ca\widehat{\boxtimes}_{\co_S} \cm
:=\limind\limits_j \ca_j\widehat{\boxtimes}_{\co_S} \cm $ on
$C_{S'}$.
Now let $\cn$ be a coherent $\co_S$-module sheaf on $S$.
Then we have that $\ca_j\widehat{\boxtimes}_{\co_S} \cn = \ca_j
\boxtimes_{\co_S} \cn$ and $\ca \widehat{\boxtimes}_{\co_S} \cn =
\ca \boxtimes_{\co_S} \cn$. Indeed, the second fact follows from the
first one, because a tensor product commute with an inductive limit.
To prove the first fact we note that it is evident when $\cn =
\co_S^{\oplus r}$, and that the functor $\ca_j
\widehat{\boxtimes}_{\co_S} (\cdot) $ is an exact functor on the
category of coherent sheaves on $S$. Now using the arguments, which
are similar to the proof of proposition 10.13 from \cite{AtM}, we
obtain that the natural map $\ca_j \boxtimes_{\co_S} \cn \lto
\ca_j\widehat{\boxtimes}_{\co_S} \cn $ is an isomorphism.
\end{rem}

Let $k$ be a field, and $\xo_{\infty} = (C, \ca)$ be a ribbon over
$k$. Let $\cb$ be a category of affine Noetherian $k$-schemes. Then
we define the following contravariant  functors
$Pic_{\sxo_{\infty}}$
 and $Pic_{X_{\infty}}$ from $\cb$ to the
category of Abelian groups.

\begin{defin} \label{picarf}
Let $\cb$ be a category of affine Noetherian $k$-schemes. Then
we define the following contravariant  functors
$Pic_{\sxo_{\infty}}$ and $Pic_{X_{\infty}}$ from $\cb$ to the
category of Abelian groups:
\begin{enumerate}
\item $ Pic_{\sxo_{\infty}}(S) \eqdef Pic(\xo_{\infty, S}) = H^1(C_S,
\ca_S^*)$;
\item $ Pic_{X_{\infty}}(S) \eqdef Pic(X_{\infty, S}) = H^1(C_S,
\ca_{S, 0}^*)$.
\end{enumerate}
\end{defin}

\subsection{Zariski tangent space.}
We recall the definition of the Zariski tangent space to a functor
at $0$. Let $\xo_{\infty} = (C, \ca)$ be a ribbon over a field $k$.
 Let $E = k \oplus k\cdot \epsilon$, where $\epsilon^2=0$, be a
$k$-algebra.
$$
T_{Pic_{\ssxo_{\infty}}}(0) \eqdef \Ker ( Pic_{\sxo_{\infty}}(\Spec
E) \lto Pic_{\sxo_{\infty}}(\Spec k)) \mbox{.}
$$
Analogously,
$$
T_{Pic_{X_{\infty}}}(0) \eqdef \Ker ( Pic_{X_{\infty}}(\Spec E) \lto
Pic_{X_{\infty}}(\Spec k)) \mbox{.}
$$

We have the following proposition.
\begin{prop}
\label{tspc} Let $\xo_{\infty} = (C, \ca)$ be a ribbon over a field
$k$.
\begin{enumerate}
\item
We have $T_{Pic_{\ssxo_{\infty}}}(0)=H^1(C,\ca )$ and
$T_{Pic_{X_{\infty}}}(0)=H^1(C,\ca_0)$.
\item
Let the ribbon $\xo_{\infty}$ correspond to some generalized
Fredholm $k$-subalgebra $A$ in $k((u))((t))$ (after a choice of a
smooth point $P \in C$ of the ribbon, formal local parameters $u,t$,
see theorem~\ref{th1}). Then we have the following exact sequence of
$k$-vector spaces:
\begin{equation} \label{tsp}
0 \lto \h^1(A) \lto T_{Pic_{X_{\infty}}}(0) \lto
T_{Pic_{\ssxo_{\infty}}}(0) \lto \h^2(A) \lto 0 \mbox{.}
\end{equation}
\end{enumerate}
\end{prop}
\Proof. Let $R = \Spec E$. We denote the base change sheaves
$$\ca'=\ca_{R}= \ca \oplus \epsilon \cdot \ca \qquad \mbox{, }\qquad \ca'_0=\ca_{R,0}=\ca
\oplus \epsilon \cdot \ca_0 \mbox{,} $$ where $\epsilon^2=0$. Then
we have canonically the following  decompositions:
$$
\ca'^*=\ca^* \times (1+ \epsilon \cdot \ca) = \ca^* \times \ca
\mbox{;}
$$
$$
\ca'^*_0=\ca_0^* \times (1+ \epsilon \cdot \ca_0) = \ca_0^* \times
\ca_0 \mbox{.}
$$
Therefore we have canonically:
$$
H^1(C, \ca'^*) = H^1(C, \ca^*) \times H^1(C, \ca) \mbox{;}
$$
$$
H^1(C, \ca'^*_0) = H^1(C, \ca^*_0) \times H^1(C, \ca_0) \mbox{.}
$$
Hence, we obtain
\begin{equation} \label{tp}
T_{Pic_{\ssxo_{\infty}}}(0)= H^1(C, \ca) \qquad \mbox{,} \qquad
T_{Pic_{X_{\infty}}}(0)= H^1(C, \ca_0) \mbox{.}
\end{equation}

We have the following exact sequence of sheaves on $C$:
$$
0 \lto \ca_0 \lto \ca \lto \ca/\ca_0 \lto 0 \mbox{.}
$$
Hence we have the following long exact sequence
\begin{equation} \label{es1}
0 \lto \frac{H^0(C, \ca/\ca_0)}{\frac{H^0(C, \ca)}{H^0(C, \ca_0)}}
\lto H^1(C, \ca_0) \lto H^1(C,\ca) \lto H^1(C, \ca/ \ca_0) \lto 0
\mbox{.}
\end{equation}
Now, using it, formulas~(\ref{tp}) and theorem~\ref{pc} we obtain
exact sequence~(\ref{tsp}).
\begin{flushright}
$\square$
\end{flushright}

\begin{rem}
\label{rem3}
According to remark~\ref{r2} and lemma~\ref{l1}, we have that if a
ribbon $\xo_{\infty}=(C, \ca)$ comes from an algebraic projective
Cohen-Macaulay surface $X$ and an ample Cartier divisor $C$ on $X$,
then exact sequence~(\ref{tsp}) transforms to the following exact
sequence:
$$
0 \lto H^1(X, \co_X) \lto T_{Pic_{X_{\infty}}}(0) \lto
T_{Pic_{\ssxo_{\infty}}}(0) \lto H^2(X, \co_X) \lto 0 \mbox{.}
$$
\end{rem}

\subsection{Formal Brauer group of an algebraic surface.}
\label{bgas}
We suppose in this subsection that a field $k$ has $\chara k =0$.

Let $X$ be a projective algebraic surface over the field $k$. We
recall the following definition of the formal Brauer group of the
surface $X$ from~\cite{AM}.
\begin{defin}
Let $\cc$ be the category of affine Artinian local $k$-schemes with
residue field $k$ (i.e. the full subcategory of affine $k$-schemes
such that $S \in \Ob (\cc)$ iff $S = \Spec B$ for an Artinian local
$k$-algebra $B$ with residue field $k$). {\em The
 formal Brauer group} $\widehat{Br}_X$ of $X$ is a contravariant functor from $\cc$
to the category of Abelian groups which is given by the following
rule:
$$
\widehat{Br}_X(S) \eqdef \Ker (H^2(X \times_k S, \co_{X \times_k
S}^*) \lto H^2(X, \co_X^*)) \mbox{,}
$$
where $S \in \Ob(\cc)$.
\end{defin}

We used the Zariski topology for the definition of the functor
$\widehat{Br}_X$. But, as it was noticed in~\cite[ch. II]{AM}
(because of the filtration with factors being coherent sheaves), we
can use, for example, the \'{e}tale topology, i.e. we have the
following equality:
$$
\widehat{Br}_X(S) \eqdef \Ker (H_{\acute{e}t}^2(X \times_k S, \co_{X
\times_k S}^*) \lto H_{\acute{e}t}^2(X, \co_X^*)) \mbox{,}
$$
where $S \in \Ob(\cc)$. It explains the name "formal Brauer group".

In~\cite[corollary 4.1]{AM} it was proved that under some conditions
on $X$ the functor $\widehat{Br}_X$ is pro-representable by the
formal group scheme ${\bf \widehat{Br}}_X$ which is a formal group
(for fields $k$ of any characteristic). It means that for any $S \in
\Ob(\cc)$:
$$
\widehat{Br}_X(S) = \Hom\nolimits_{form. sch.} (S, {\bf
\widehat{Br}}_X) \mbox{,}
$$
where $\Hom_{form. sch.}$  is considered in the category of formal
schemes.

Since we supposed that $\chara k =0$, we will give an easy proof
that the functor $\widehat{Br}_X$ is always pro-represantable in the
following lemma.
\begin{lemma}
The functor $\widehat{Br}_X$ from the category $\cc$ to the category
of Abelian groups is pro-representable by the formal group scheme
${\bf \widehat{Br}}_X = \Spf \widehat{Sym}_k(H^2(X, \co_X)^*)$,
where the group law in the formal group ${\bf \widehat{Br}}_X$ is
given by $v \longmapsto v \otimes 1 + 1 \otimes v$, $v \in H^2(X,
\co_X)^*$.
\end{lemma}
\Proof. By definition, we have for $k$-algebra
$$
\widehat{Sym}_k (H^2(X, \co_X)^*) \eqdef \prod_{i =0}^{\infty}
S^i(H^2(X, \co_X)^*) \mbox{,}
$$
where $S^i(\cdot)$ is the $k$-vector space of symmetric $i$-th
tensors over the field $k$, $S^0(\cdot) = k$. $ \widehat{Sym}_k
(H^2(X, \co_X)^*)$ is a topological local $k$-algebra over a
discrete field $k$. This topology is given by the infinite product
topology of discrete spaces.

For any $S = \Spec B \in \Ob(\cc)$ we have $B = k \oplus I$, where
$I$ is the maximal ideal in the ring $B$, $\dim_k I < \infty$ and
$I^n =0$ for some $n \ge 0$. We consider the discrete topology on
the ring $B$. Therefore we have
$$
\Hom\nolimits_{form. sch.} (S, \Spf \widehat{Sym}_k(H^2(X,
\co_X)^*))= \Hom\nolimits_{k-alg., cont} (\widehat{Sym}_k(H^2(X,
\co_X)^*), B) =
$$
$$
= \Hom\nolimits_k (H^2(X, \co_X)^*, I )= H^2(X, \co_X) \otimes_k I
\mbox{,}
$$
where $\Hom_{k-alg., cont.}$ is considered in the category of
topological $k$-algebras. We used also that ${H^2(X, \co_X)^*}^*=
H^2(X, \co_X)$, because $\dim_k H^2(X, \co_X) < \infty$.

On the other hand, for base change sheaf $\co'_X=\co_{X \times_k S}
$ on $X$ we have
\newline
$\co'_X= \co_X \oplus (\co_X \otimes_k I)$. Hence $\co'^*_X =
\co^*_X \times (1 + \co_X \otimes_k I)$. The exponential map gives
an isomorphism of sheaves of Abelian groups:
$$
exp \quad : \quad \co_X \otimes_k I \lto 1 + \co_X \otimes_k I
\mbox{.}
$$
Therefore $\co'^*_X= \co^*_X \times (\co_X \otimes_k I)$. Hence
$$
H^2(X, \co'^*_X)= H^2(X, \co^*_X) \times H^2(X, \co_X \otimes_k I) =
H^2(X, \co^*_X) \times (H^2(X, \co_X) \otimes_k I) \mbox{.}
$$
Therefore we have
$$
\widehat{Br}_X(S) = H^2(X, \co_X) \otimes_k I \mbox{.}
$$
\begin{flushright}
$\square$
\end{flushright}

\vspace{0.5cm}

Moreover, we have the universal object
$$
\tau \; \in \; \Ker(H^2(X \times_k {\bf \widehat{Br}}_X, \co^*_{X
\times_k {\bf \widehat{Br}}_X}) \to H^2(X, \co_X^*)) \mbox{,}
$$
which is constructed in the following way. We have
$$
\Ker(H^2(X \times_k {\bf \widehat{Br}}_X, \co^*_{X \times_k {\bf
\widehat{Br}}_X}) \to H^2(X, \co_X^*)) = H^2(X, 1+J
\mathbin{\hat{\otimes}}_k \co_X ) = H^2(X, \co_X) \otimes_k J
\mbox{,}
$$
where the ideal $J = \prod\limits_{i=1}^{\infty} S^i(H^2(X,
\co_X^*))$. Besides,  the sheaf of Abelian groups $J
\mathbin{\hat{\otimes}}_k \co_X $ is isomorphic to the sheaf $1 + J
\mathbin{\hat{\otimes}}_k \co_X $ via the exponential map.

Now $\tau =  Id$, where ${Id}$ is the identity map from
$$\End\nolimits_k(H^2(X, \co_X))= H^2(X, \co_X) \otimes_k H^2(X, \co_X)^* \mbox{.}$$
And there is a canonical embedding of $k$-vector spaces:
$$
H^2(X, \co_X) \otimes_k H^2(X, \co_X)^* \subset H^2(X, \co_X)
\otimes_k J \mbox{.}
$$

\subsection{Formal Brauer group of a ribbon.}
\label{bgr}
We use in this subsection the same notations as in
subsection~\ref{bgas}. In particularly, a field $k$ has $\chara
k=0$, $\cc$ is the category of affine Artinian local $k$-schemes
with residue field $k$. We introduce the following definition.

\begin{defin} Let $\xo_{\infty}=(C, \ca)$ be a ribbon over  a field
$k$. {\em The
 formal Brauer group} $\widehat{Br}_{\sxo_{\infty}}$ of $\xo_{\infty}$ is a contravariant functor from $\cc$
to the category of Abelian groups which is given by the following
rule:
$$
\widehat{Br}_{\sxo_{\infty}}(S) \eqdef \Ker (H^1(C_S,
\ca_S^*/\ca_{S, 0}^*) \lto H^1(C, \ca^*/ \ca_0^*)) \mbox{,}
$$
where $S \in \Ob(\cc)$.
\end{defin}

Below, in remark~\ref{rmr},  we will explain, why we use the name
"formal Brauer group" for a ribbon. Now we have the following
proposition.
\begin{prop}
\label{prop3} Assume $C$ is a projective curve. Then the functor
$\widehat{Br}_{\sxo_{\infty}}$ from the category $\cc$ to the
category of Abelian groups is pro-representable by a formal group
scheme ${\bf \widehat{Br}_{\sxo_{\infty}}}$.
\end{prop}
\Proof. We denote a $k$-vector space $V = H^1(C, \ca/ \ca_0) $. (We
note that, by theorem~\ref{pc}, $V = \h^2(A)$ when $A$ is a
generalized Fredholm subspace in $k((u))((t))$ which correspond to
the ribbon $\xo_{\infty}=(C, \ca)$  with some smooth point $P \in C$
of the ribbon, formal local parameters $u,t$, see
theorem~\ref{th1}). We note that we have canonically:
$$
V = \limind_{i \ge 0} V_i \mbox{,}
$$
where $V_i = H^1(C, \ca_{-i} / \ca_0)$, $i \ge 0$. And $\dim_k V_i <
\infty $, since $C$ is a projective curve and $\ca_{-i} / \ca_0$
are coherent sheaves on the scheme $X_{i-1}$ by \cite[prop. 1]{Ku}.
Therefore we have
$$
V^* = \limproj_{i \ge 0} V_i^*  \mbox{.}
$$
The $k$-vector space $V^*$ has a natural linearly compact topology,
which is given by topology of this projective limit, where every
$V_i^*$ has a discrete topology, see \cite[ch.III, \S 2, ex. 15]{Bu}.

For any $l \ge 0$ we define
$$
S_{cont}^l(V^*) \eqdef  \limproj_{i \ge 0} S^l(V_i^*) \mbox{.}
$$
These spaces has also a linearly compact topology, which is  given
by the projective limit. We define
$$
T = \widehat{Sym}_{k, cont}(V^*)  \eqdef  \prod_{l=0}^{\infty}
S_{cont}^l (V^*) \mbox{.}
$$
By construction, we have $S_{cont}^0 (V^*) =k$. The $k$-vector space
$T$ has the product topology. Therefore $T$ is a linearly compact
space as the product of linearly compact spaces. Hence $T$ is Hausdorff and
complete, see \cite[ch.III, \S 2, ex. 16]{Bu}.

For any $l_1 \ge 0$, $l_2 \ge 0$ we have canonical continuous
bilinear  map over $k$:
$$
S_{cont}^{l_1}(V^*) \times S_{cont}^{l_2}(V^*) \lto
S_{cont}^{l_1+l_2}(V^*) \mbox{.}
$$
Therefore $T$ is a topological local $k$-algebra (over  the discrete
field $k$). The maximal ideal $J$ in $T$ is given as
$$
J \eqdef \prod_{l=1}^{\infty} S_{cont}^l (V^*) \mbox{.}
$$
By construction, for any open $k$-subspace $U \subset T$ there is $m
>0$ such that $U \supset J^m$. Using these properties of topological $k$-algebra $T$
we obtain that the following formal scheme is well-defined
(see~\cite[ch. I, \S 10]{EGAI}):
$$
{\bf \widehat{Br}_{\sxo_{\infty}}} \eqdef \Spf (T) \mbox{.}
$$
Moreover, ${\bf \widehat{Br}_{\sxo_{\infty}}}$ is a formal group
with the group law $v \mapsto v \otimes 1 + 1 \otimes v$,
\newline
$v \in V^*= S^1_{cont}(V^*)$. ($V^*$  topologically generates the
$k$-algebra $T$.)

Now we have to check that the formal group scheme ${\bf
\widehat{Br}_{\sxo_{\infty}}}$ pro-represents the functor
$\widehat{Br}_{\sxo_{\infty}}$.

For any $S = \Spec B \in \Ob(\cc)$ we have $B = k \oplus I$, where
$I$ is the maximal ideal in the ring $B$, $\dim_k I < \infty$ and
$I^n =0$ for some $n \ge 0$. We consider the discrete topology on
the ring $B$. Therefore  we have
$$
\Hom\nolimits_{form. sch.} (S, {\bf \widehat{Br}_{\sxo_{\infty}}})=
\Hom\nolimits_{k-alg, cont} (T, B)  = \Hom\nolimits_{k, cont} (V^*,
I ) \mbox{,}
$$
where $\Hom\nolimits_{k, cont}$ is considered in the category of
topological $k$-vector spaces. Since $\dim_k(V_i) < \infty$, $i \ge
0$, we have
\begin{equation} \label{eq5}
\Hom\nolimits_{k, cont} (V^*, k ) = V \mbox{.}
\end{equation}
(We note also that $V^* = \Hom\nolimits_{k, cont} (V, k )$, where
$V$ has a discrete topology.) Since $\dim_k I < \infty$ and $I$ has
a discrete topology, we obtain the following formula from
formula~(\ref{eq5}):
$$
\Hom\nolimits_{k, cont} (V^*, I ) =V \otimes_k I \mbox{.}
$$
Therefore
\begin{equation} \label{fb}
 \Hom\nolimits_{form. sch.} (S, {\bf
\widehat{Br}_{\sxo_{\infty}}})=V \otimes_k I \mbox{.}
\end{equation}

On the other hand, we have the following split exact sequence:
$$
1\rightarrow 1+I\otimes_k \ca /\ca_0\rightarrow \ca^*_S/\ca^*_{S,0} \rightarrow \ca^* /\ca_0^*\rightarrow 0,
$$
which is the factor of the exact sequence
\begin{equation}
\label{ca} 1\rightarrow 1+I\otimes_k \ca \rightarrow \ca^*_S
\rightarrow \ca^* \rightarrow 0
\end{equation}
by the following exact sequence
\begin{equation}
\label{ca_0}
1\rightarrow 1+I\otimes_k \ca_0 \rightarrow \ca^*_{S,0} \rightarrow \ca^*_0 \rightarrow 0.
\end{equation}
The sheaf of Abelian groups $I\otimes_k \ca /\ca_0$ is isomorphic to
the sheaf $1+I\otimes_k \ca /\ca_0$ via the exponential map.
Therefore, we have
$$
\widehat{Br}_{\sxo_{\infty}}(S)= H^1(C, 1+I\otimes_k \ca
/\ca_0)\simeq H^1(C, I\otimes_k \ca /\ca_0)=I\otimes_kH^1(C
\otimes_k \ca /\ca_0) \mbox{.}
$$
So, $\widehat{Br}_{\sxo_{\infty}}(S)= \Hom\nolimits_{form. sch.} (S,
{\bf \widehat{Br}_{\sxo_{\infty}}}) $.

\begin{flushright}
$\square$
\end{flushright}

\begin{rem}
\label{fb1}
Let $A$ be any commutative $k$-algebra. Then we have an analog of
formula~(\ref{fb}):
$$
 \Hom\nolimits_{form. sch.} (\Spec A, {\bf
\widehat{Br}_{\sxo_{\infty}}})=V \otimes_k N_A  \mbox{,}
$$
where $N_A$ is the nilradical of the ring $A$. Indeed, following the
proof of formula~(\ref{fb}), we see that it is enough to prove the
following formula
\begin{equation} \label{eqq5}
\Hom\nolimits_{k, cont} (V^*, N_A ) =V \otimes_k N_A \mbox{,}
\end{equation}
where $N_A$ has a discrete topology. But we have that $V^*$ is a
linearly compact $k$-vector space. Therefore for any $\phi \in
\Hom\nolimits_{k, cont} (V^*, N_A )$ we have that $\phi(V^*)$ is a
linearly compact $k$-vector subspace in a discrete $k$-vector space
$N_A$. Hence $\dim_k \phi(V^*) < \infty$. Now, using
formula~(\ref{eq5}), we obtain formula~(\ref{eqq5}).
\end{rem}

\begin{rem} \label{rmr}
By remark~\ref{r2}  and theorem~\ref{pc} we have that if a ribbon
$\xo_{\infty}=(C, \ca)$ comes from an algebraic projective
Cohen-Macaulay surface $X$ and an ample Cartier divisor $C$ on $X$,
then $H^2(X, \co_X) = H^1(C, \ca/ \ca_0)$. Therefore we have
$$
{\bf \widehat{Br}}_X = \Spf \widehat{Sym}_k(H^2(X, \co_X)^*) = {\bf
\widehat{Br}_{\sxo_{\infty}}} \mbox{.}
$$
\end{rem}

\subsection{Formal Picard group of a ribbon.}

We use in this subsection the same notations as in
subsection~\ref{bgr}.

\begin{defin} Let $\xo_{\infty}=(C, \ca)$ be a ribbon over  a field
$k$. {\em The
 formal Picard group} $\widehat{Pic}_{\sxo_{\infty}}$ of $\xo_{\infty}$ is a contravariant functor from $\cc$
to the category of Abelian groups which is given by the following
rule:
$$
\widehat{Pic}_{\sxo_{\infty}}(S) \eqdef \Ker (H^1(C_S,
\ca_S^*) \lto H^1(C, \ca^*)) \mbox{,}
$$
where $S \in \Ob(\cc)$.

Analogously, {\em the formal Picard group}
$\widehat{Pic}_{X_{\infty}}$ of $X_{\infty}$ is a functor from $\cc$
to the category of Abelian groups given by
$$
\widehat{Pic}_{X_{\infty}}(S) \eqdef \Ker (H^1(C_S,
\ca_{S,0}^*) \lto H^1(C, \ca^*_0)).
$$
\end{defin}

Let $S = \Spec B \in \Ob (\cc)$, where $B=k \oplus I$, $I$ is the
maximal ideal of the ring $B$, $\dim_k I < \infty$, and $I^n =0$ for
some $n \ge 0$.

As it follows from  split exact sequences (\ref{ca}), (\ref{ca_0}),
using exponential and logarithmic maps,  we have
$$\widehat{Pic}_{\sxo_{\infty}}(S)=I\otimes_k
H^1(C,\ca ) \mbox{,} \qquad \widehat{Pic}_{X_{\infty}}(S)=I\otimes_k
H^1(C,\ca_0 ) \mbox{.}$$ So, if we define the contravariant functor
$P$ from $\cc$ to the category of Abelian groups by the rule
$$
P(S)\eqdef I\otimes_k H, \mbox{\quad where\quad } H\eqdef\frac{H^0(C,
\ca/\ca_0)}{\frac{H^0(C, \ca)}{H^0(C, \ca_0)}} \mbox{,}
$$
we get the following exact sequence of groups, which is functorial
on $\cc$ (compare with sequence (\ref{es1})):
\begin{equation}
\label{fseq}
0\lto P(S)\lto \widehat{Pic}_{X_{\infty}}(S) \lto \widehat{Pic}_{\sxo_{\infty}}(S)
\lto \widehat{Br}_{\sxo_{\infty}}(S)\lto 0
\end{equation}
We define another functor $\overline{Pic}$ from $\cc$ to the
category of Abelian groups by the rule
$$
\overline{Pic}(S)\eqdef I\otimes_k \frac{H^1(C,\ca_0 )}{H} \mbox{.}
$$
(If $H =0$, then $\overline{Pic}= \widehat{Pic}_{X_{\infty}} $.)
Then from  sequence (\ref{fseq}) we obtain another exact sequence:
\begin{equation}
\label{fseq1} 0\lto \overline{Pic}(S)\lto
\widehat{Pic}_{\sxo_{\infty}}(S)\lto
\widehat{Br}_{\sxo_{\infty}}(S)\lto 0
\end{equation}

\begin{prop}
\label{prop4}
Assume $C$ is a projective curve. Then we have:
\begin{enumerate}
\item
\label{1)} There is a noncanonical functorial (with $S$) splitting
of sequence~(\ref{fseq1}).
\item
\label{2)} If $\dim_k H<\infty$, then the functor $\overline{Pic}$
from the category $\cc$ to the category of Abelian groups is
pro-representable by the formal group scheme ${\bf \overline{Pic}}$.
\end{enumerate}
\end{prop}
\begin{rem}
The condition $\dim_k H<\infty$ of this proposition is satisfied,
for example, if the ribbon comes from an algebraic projective
surface $X$ as in remark \ref{rem3}, since in this case
$H=H^1(X,\co_X)$. Another example is a ribbon coming from the Schur
pair $(A,W)$ as in theorem \ref{th1}, where $A$ is chosen so that
$\dim_k \h^1(A)<\infty$. Due to lemma \ref{l1} one can easily
construct a lot of examples of such spaces.
\end{rem}

\Proof. The first claim is clear, because we can fix any $k$-linear
section of the map $ H^1(C,\ca ) \rightarrow H^1(C,\ca /\ca_0)$ and
then extend it for any $S$ in~(\ref{fseq1}) by tensor product with
identity map on $I$ over $k$.

The proof of the second assertion is similar to the proof of proposition \ref{prop3}.
Namely, let $\dim_k H<\infty$. We have $H^1(C,\ca_0 )= \limproj\limits_{j>0} H^1(C,\ca_0/\ca_{j})$
by \cite[corollary 1]{Ku}, and $\dim_k H^1(C,\ca_0/\ca_{j})<\infty$.
The homomorphism $i:H\rightarrow H^1(C,\ca_0 )$ gives a system of compatible
homomorphisms $i_j:H\rightarrow H^1(C,\ca_0/\ca_{j})$.
Denote by $K_j$ the kernel of $i_j$ and by $C_j$ the cokernel of $i_j$.
Then we obtain the exact sequence of projective systems of $k$-vector spaces:
$$
0\lto (K_j)\lto (H_j)\lto H^1(C,\ca_0/\ca_{j}) \lto (C_j)\lto 0,
$$
where $H_j=H$ for all $j$. Since $\dim_k K_j<\infty$ for all $j$,
the systems $(K_j)$, $(H_j)$ satisfy the Mittag-Leffler condition.
Then by \cite[lemma 1]{Ku} we obtain that
$$H^1(C,\ca_0 )/H\simeq \limproj\limits_{j\in \sdn}C_j \mbox{,}$$
where $\dim_k C_j<\infty$ for all $j$. Denote by $V$ the $k$-vector
space $H^1(C,\ca_0 )/H$.

We have $V_{cont}^*:= \Hom_{k,cont}(V, k)=\limind\limits_{j}C_j^*$
is a $k$-vector space with a discrete topology. Now we define
$$
T:=\widehat{Sym}_k(V_{cont}^*)=\prod_{l=0}^{\infty}S^l(V_{cont}^*).
$$
$T$ is a topological local $k$-algebra with the maximal ideal
$J=\prod\limits_{l=1}^{\infty}S^l(V_{cont}^*)$. The topology on $T$
is a linear product topology. It is clear that $J$ is the maximal
ideal of definition and that $T$ is an admissible ring (and
moreover, adic) in the sense of \cite[0.7.1]{EGAI}. Therefore we can
define
$$
{\bf \overline{Pic}}\eqdef \Spf (T).
$$
Again, as in proposition \ref{prop3}, $ {\bf \overline{Pic}}$ is a
formal group with the group law $v \longmapsto v\otimes 1 + 1
\otimes v$, $v \in V_{cont}^*$.

For any $S = \Spec B \in \Ob(\cc)$ we have
$$
\Hom\nolimits_{form. sch.} (S, {\bf \overline{Pic}})=
\Hom\nolimits_{k-alg, cont} (T, B)  = \Hom\nolimits_{k, cont}
(V_{cont}^*, I )= V\otimes_k I=\overline{Pic}(S).
$$

\begin{flushright}
$\square$
\end{flushright}

Now we have the following obvious corollary of
proposition~\ref{prop4}.
\begin{corol}
\label{cor1}
Let $C$ be a projective curve, and $\dim_k H<\infty$. Then
the functor $\widehat{Pic}_{\sxo_{\infty}}$ is
(noncanonically) pro-representable by the formal group scheme ${\bf
\overline{Pic}}\times_k{\bf \widehat{Br}_{\sxo_{\infty}}}$. Such
decompositions are in one-to-one correspondence with functorial
(with $S$) splittings of sequence~(\ref{fseq1}).
\end{corol}

\section{Picard functor of a ribbon.}
\label{picfunrib}
\subsection{The order function}
\label{funorder}

In this subsection we will give an appropriate generalization of the order function used in \cite[\S 4]{Ku}. It will be used later.

For a topological space $U$ let $W_U(\dz )$ be a sheaf of functions on $U$ with values in $\dz$.

Let $\xo_{\infty}$ be a ribbon over a Noetherian scheme $S$.

{\it In the following we assume that for every} $s\in S$ {\it there exists a point}
$P_s\in \xo_{\infty, s}$ {\it such that} $(\ca_{s,1})_{P_s}(\ca_{s,-1})_{P_s}=(\ca_{s,0})_{P_s}$ {\it and that the underlying topological
space of} $\xo_{\infty, s}$ {\it is irreducible}.
{\it We also assume that the morphism $\tau :C\rightarrow S$ from definition~\ref{ribbon} is  locally of finite type}.

Note that by \cite[prop.9]{Ku} the function of order ord defined in \cite[def.6]{Ku} is a morphism of sheaves of
groups on any $\xo_{\infty, s}$.

\begin{rem}
\label{remorder} For example, this assumption is satisfied in the
case of the ribbon $\xo_{\infty ,S}$, where $S \rightarrow \Spec k$
is a base change, and $ \xo_{\infty}$ is a ribbon over an
algebraically closed field $k$ with irreducible underlying
topological space and either with a smooth point in the sense of
\cite[def. 9]{Ku}, or with condition $(**)$ from
definition~\ref{uslovie}.

Indeed, in this case for every $s\in S$ we have that the underlying topological space of $\xo_{\infty, s}$ is an
irreducible curve by \cite[vol.I, ch.III, \S 15, th.40, cor.1]{ZS} (see also \cite[ch. II, ex.3.20]{Ha}),
and $\tau$ is of finite type. If $P$ is a smooth point of the ribbon $ \xo_{\infty}$, and $P_s$ is the closed
point that maps to $P$, then $P_s$ is a smooth point of the ribbon $\xo_{\infty, s}$. The reason is that we
can lift the elements $t\in \ca_{1,P}$, $t'\in \ca_{-1,P}$ with $tt'=1$ to analogous
elements $t_s\in \ca_{1,P_s}$, $t_s'\in \ca_{-1,P_s}$. Then, for example, the arguments from the
proof of \cite[prop.7]{Ku} show that $P_s$ is smooth.
The same arguments work in the case of condition $(**)$.
\end{rem}

\begin{defin}[order map]
\label{order}
Define a morphism of sheaves of sets
$$
 \ord :\ca^*\lto W_C(\dz ), \mbox{\quad } \ord (a)(x)\eqdef \max\{j| \mbox{\quad } a|_{U_s}\in \ca_{s,j}(U_s)\},
$$
where $a\in \ca^*(U)$ for an open $U\subset C$, $x\in U$, $s=\tau (x)$.

On any $\xo_{\infty, s}$ $\ord$ coincide with  the order function from \cite[def. 6]{Ku}.
By \cite[prop.~9]{Ku} $\ord$ is compatible with restriction homomorphisms on any $\xo_{\infty, s}$. So, our definition is correct.
\end{defin}

We want to give a condition when $\ord$ is a morphism of sheaves of groups, and when it factors through the sheaf $\dz_C\subset W_C(\dz )$ of locally constant functions. We also want to describe in this case the kernel of the order map.

If $S=\Spec K$, where $K$ is a field, then this definition coincide with the definition 6 of \cite{Ku}. In \cite[prop.~8, prop.~9]{Ku} we gave certain sufficient conditions for the order function to be a homomorphism (obviously, in this case it is locally constant), see also counter-example 7 in {\it loc.cit.}.

\begin{lemma}
\label{tex} In our assumptions we have: for every $P_s$ there exists
its affine neigbourhood $U_{P_s} \subset C$ such that all
$\ca_j|_{U_{P_s}}$ are invertible sheaves of
$\ca_0|_{U_{P_s}}$-modules and
$\ca_{-j}|_{U_{P_s}}=\ca_j^{-1}|_{U_{P_s}}$.
\end{lemma}

\Proof. First, let's prove that the natural homomorphism of $\co_{C,P_s}$-modules
\begin{equation}
\label{qwert}
(\ca_{-1}/\ca_0)_{P_s}\otimes_{\co_{C,P_s}}(\ca_1/\ca_2)_{P_s}\lto (\ca_0/\ca_1)_{P_s}=\co_{C, P_s}
\end{equation}
is an isomorphism.

Since by our assumption there is an isomorphism
$$[(\ca_{-1}/\ca_0)_{P_s}\otimes_{\co_{S,s}} k(s)]\otimes_{\co_{C,s}}[(\ca_1/\ca_2)_{P_s}\otimes_{\co_{S,s}} k(s)]\simeq \co_{C_s,P_s},
$$
the homomorphism (\ref{qwert}) is surjective by Nakayama's lemma. Let $K$ be the kernel of this map. Then, since $\co_{C, P_s}$ is a flat $\co_{S,s}$-module, the following sequence is exact (by \cite[ch.2, ex.26]{AtM}):
$$
0\lto K\otimes_{\co_{S,s}}k(s)\lto ((\ca_{-1}/\ca_0)_{P_s}\otimes_{\co_{C,P_s}}(\ca_1/\ca_2)_{P_s})\otimes_{\co_{S,s}}k(s)\lto \co_{C_s,P_s}\lto 0.
$$
Since
$$((\ca_{-1}/\ca_0)_{P_s}\otimes_{\co_{C,P_s}}(\ca_1/\ca_2)_{P_s})\otimes_{\co_{S,s}}k(s)\simeq [(\ca_{-1}/\ca_0)_{P_s}\otimes_{\co_{S,s}} k(s)]\otimes_{\co_{C,s}}[(\ca_1/\ca_2)_{P_s}\otimes_{\co_{S,s}} k(s)],
$$
we get $0=K\otimes_{\co_{S,s}}k(s)=K/\cm_sK$. Therefore, $K=0$ by
Nakayama's lemma.

Now let $U_{P_s}$ be an affine neigbourhood, where there exist
$\bar{t}'\in  \ca_{-1}/\ca_0(U_{P_s})$, $\bar{t}\in
\ca_1/\ca_2(U_{P_s})$ such that $\bar{t}\bar{t}'=1$. Then, by
\cite[prop.3]{Ku}, we can lift the elements $\bar{t}, \bar{t}'$ and
find $t'\in \ca_{-1}(U_{P_s})$, $t\in \ca_1(U_{P_s})$ such that
$tt'= 1$. Then for any $j$ we have
$\ca_j|_{U_{P_s}}=t^j(\ca_0|_{U_{P_s}})$ (compare with the arguments
in the proof of prop.~7 from~\cite{Ku}).
\begin{flushright}
$\square$
\end{flushright}

\begin{rem}
 \label{osobu} If for any $s\in S$ the ribbon $\xo_{\infty ,s}$
satisfies the condition $(**)$ from definition \ref{uslovie}, then
the statements of lemma are valid for the sheaves $\ca_j$ on the
whole space $C$ (not only on $U_{P_s}$). The proof is the same.
\end{rem}

Let's consider several cases.
\vspace{0.3cm}

\noindent {\bf Case 1.} Let $S$ be an integral scheme. We claim that the order map on $\ca^*|_{U_{P_s}}$ factors through $\dz_C|_{U_{P_s}}$ and is a morphism of sheaves of Abelian groups. Moreover, $(\ca^*/\ca_0^*)|_{U_{P_s}}\simeq \dz_C|_{U_{P_s}}$.

Let $a\in \ca^*(U_{P_s})$ and $j$ be the biggest integer with $a\in \ca_j(U_{P_s})$. Then $j=\ord (a)(x)$ for any $x\in U_{P_s}$. Indeed, by lemma \ref{tex} there exists an invertible element $t\in \ca_1(U_{P_s})$. So, $a=a_0t^j$ with $a_0\in \ca_0(U_{P_s})\backslash \ca_1(U_{P_s})$. If $a^{-1}\in \ca_k(U_{P_s})\backslash \ca_{k+1}(U_{P_s})$, then $a^{-1}=b_0t^k$, $b_0\in \ca_0(U_{P_s})\backslash \ca_1(U_{P_s})$. Then $1=a^{-1}a=a_0b_0t^{j+k}$, hence $j+k\le 0$ and $a_0b_0=t^{-j-k}\notin \ca_1(U_{P_s})$, because $\ca_0(U_{P_s})/ \ca_1(U_{P_s})\simeq \co_C(U_{P_s})$ has no zero divisors, since $C$ is irreducible and reduced, what follows from our assumptions (see remark \ref{rem10} below).

Therefore, $j+k=0$, $b_0=a_0^{-1}$ and $a_0\in \ca_0^*(U_{P_s})$, $a=a_0t^j$, $b^{-1}=a_0^{-1}t^{-j}$. Clearly, this is preserved under base change $s\rightarrow S$, therefore $j=\ord (a)(x)$ for any $x\in U_{P_s}$. So,  the degree map factors through $\dz_C|_{U_{P_s}}$ and is, obviously, a morphism of sheaves of Abelian groups with $(\ca^*/\ca_0^*)|_{U_{P_s}}\simeq \dz_C|_{U_{P_s}}$, since $\ord (t)|_{U_{P_s}}=1$.
\begin{rem}
\label{rem10}
$C$ is irreducible, because $S$ is irreducible. Indeed, assume the converse. Then there are two open subsets $U_1\subset C$, $U_2\subset C$ with $U_1\cap U_2=\varnothing $. Since $\tau :C\rightarrow S$ is flat and locally of finite type, it is open and therefore $\tau (U_1)\cap \tau (U_2)\neq \varnothing$. So, if $s\in \tau (U_1)\cap \tau (U_2)$, then $C_s\cap U_1\neq \varnothing$, $C_s\cap U_2\neq \varnothing$ and therefore $C_s$ is reducible, a contradiction with our assumption.

To prove that $C$ is reduced, let's assume the converse. We can assume $S$ is affine and the nilradicals $Nil(\co_C(U))\ne 0$ for any open $U$. Let $S'$ be a normalization of $S$ and $\tau' :C'\rightarrow S'$ be the base change. Since $C$ is flat over $S$, we
have $Nil(\co_{C'}(U\times_S S'))\ne 0$ for any affine $U\subset C$, because we have the embedding $Nil(\co_C(U))\hookrightarrow Nil(\co_{C'}(U\times_SS'))$. For any point $s\in S'$ of codimension 1 let $T\in C'$ with $\tau '(T)=s$. Then we have $Nil(\co_{C',T})\ne 0$. But $\co_{C',T}$ is a flat $\co_{S',s}$-module and $\co_{S',s}$ is a regular local ring of dimension 1. Moreover, $\co_{C',T}\otimes k(s)\simeq \co_{C',T}/u\co_{C',T}\simeq \co_{C'_s,T}$, where $u$ is a generator of the maximal ideal of $\co_{S',s}$, has no zero divisor, because by our assumptions $C_s$ is irreducible curve. Therefore, $\co_{C',T}/u\co_{C',T}\simeq \co_{C',T,red}/u\co_{C',T,red}$, where $\co_{C',T,red}=\co_{C',T}/Nil(\co_{C',T})$. Note also that $u$ is not a zero divisor in $\co_{C',T,red}$, since it is not a zero divisor in $\co_{C',T}$ by the local flatness criterium (\cite[ch.III, \S 5, th.1]{Bu} or \cite[ch.III, lemma 10.3.A]{Ha}). Therefore, $\co_{C',T,red}$ is a flat $\co_{S',s}$-module by this criterium. Hence,  $Nil(\co_{C',T})\otimes k(s)=0$ and by the Nakayama lemma $Nil(\co_{C',T})=0$, a contradiction.
\end{rem}

\begin{rem}
In situation of remark \ref{osobu}
the statements of our case are valid for the whole space $C$ (not only on $U_{P_s}$). The proof is the same.
\end{rem}

Now we claim that the order map on $\ca^*$ factors through $\dz_C$ on the whole space $C$ (although may be $(\ca^*/\ca_0^*)\not\simeq \dz_C$).

Indeed, let $U$ be a neigbourhood of a point $x\in C$, $a\in \ca^*(U)$.  Then by \cite[prop.9]{Ku} and by definition, for all points $y\in U_s$, where $s=\tau (x)$, we have $\ord (a)(y)=\ord (a)(x)$, because $C_s$ is irreducible. Since $C$ is irreducible, we have $U\cap U_{P_q}\neq \varnothing$ for any $q\in S$. Then, by the arguments above, we have $\ord (a)(x)= \ord (a)(y)$ for  any $y\in U\cap U_{P_s}$. Analogously, for $x'\in U$, $s'=\tau (x')$ we have $\ord (a)(x')= \ord (a)(y)$ for  any $y\in U\cap U_{P_{s'}}$. Since $U\cap U_{P_{s'}}\cap U_{P_s}\neq \varnothing$, we obtain $\ord (a)(x)=\ord (a)(x')$ for any $x'\in U$.
\begin{rem}
\label{ordkernel}
We can define in our case the map $ord$ as in \cite{Ku}:
$$
ord (a)\eqdef \max_{j\in \sdz} \{ j| \mbox{\quad} a\in \ca_j(U)\},
$$
where $a\in \ca^*(U)$. Then we claim that $ord (a)= \ord (a)(x)$ for any $x\in U$.

Indeed, we have proved above that $\ord (a)(x)=\ord (a)(x')$ for any $x'\in U$ and
$$
\ord (a|_{U_{P_s}\cap U})(y)= ord (a|_{U_{P_s}\cap U})\ge ord (a)
$$
for any $y\in U_{P_s}\cap U$. If $ord (a|_{U_{P_s}\cap U})> ord (a)$, then this would mean that that the image of the element $\bar{a}\in \ca_{ord (a)}/\ca_{ord (a)+1}(U)$ under the map $\varphi : \ca_{ord (a)}/\ca_{ord (a)+1}(U) \rightarrow (\ca_{ord (a)}/\ca_{ord (a)+1})_{\eta}(U_{\eta})$, where $\eta$ is a general point on $S$, is zero. But $\varphi$ is an injective map, because $\ca_{ord (a)}/\ca_{ord (a)+1}(U)$ is a flat $\co_S(\tau (U))$-module, $ \ca_{ord (a)}/\ca_{ord (a)+1}$ is a coherent sheaf, and the map $\co_S(\tau (U))\rightarrow \co_{S,\eta }$ is an embedding. So, $ord (a|_{U_{P_s}\cap U})= ord (a)$ and $ord (a)= \ord (a)(x)$ for any $x\in U$.
\end{rem}

\noindent
{\bf Case 2.} Let $S$ be a reduced scheme. We claim that the same assertions as in Case 1 hold.

Let $a\in \ca^*(U_{P_s})$ and $j=\ord (a)(x)$, where we can assume $x\in C$ to be a point such that $s=\tau (x)$ belong to several irreducible components $S_1,\ldots S_k$ (without loss of generality we can assume $S=S_1\cup\ldots\cup S_k$). By Case 1 we know that $a_0=at^{-j}|_{C_{S_i}}\in \ca^*_{ S_i,0}(U_{P_s,S_i})$. For any $l,m\in \dz$ we have the exact sequences of sheaves of filtered $\ca_0$-modules
$$
0\lto \ca_{m,S_1\cup (S_2\cup\ldots \cup S_l)}\lto \ca_{m,S_1}\times\ca_{m,(S_2\cup\ldots \cup S_l)}\lto \ca_{m,S_1\cap (S_2\cup\ldots \cup S_l)}\lto 0.
$$
Therefore, by obvious induction arguments, using the exact sequences, we obtain  $a_0\in \ca_0(U_{P_s})$. Similarly, $b_0=a^{-1}t^j\in \ca_0(U_{P_s})$ and $a_0b_0=1$. So, $a_0,b_0\in \ca_0^*(U_{P_s})$.

Thus $\ord (a)$ is locally constant on $U_{P_s}$.

\begin{rem}
In situation of remark \ref{osobu}
the statements of our case are valid for the whole space $C$ (not only on $U_{P_s}$). The proof is the same.
\end{rem}

To show that $\ord (a)$ is locally constant on any $U\subset C$, we can repeat the arguments from the end of Case 1, because $C_s$ is irreducible (so, $U\cap U_{P_s}\neq \varnothing$, and $U\cap U_{P_s}\cap S_i\times_S C\neq \varnothing$ for any $i=1,\ldots k$ (so, $U\cap U_{P_{s'}}\cap U_{P_s}\neq \varnothing$ for $s'\in S_i$, $k\in \{1,\ldots ,k\}$)).

Note that, using remark \ref{ordkernel} and above arguments, we obtain that an element $a\in \ca^*(U)$ with $\ord (a)\equiv 0$ must belong to $\ca_0^*(U)$.
\vspace{0.3cm}

\noindent {\bf Case 3.} Let $S$ be an arbitrary Noetherian scheme.
Let $\cn_S\eqdef \Ker (\co_S\rightarrow \co_{S_{red}})$ be the
nilradical. Since $\cn_S$ is a coherent sheaf on $S$, we have by the
arguments of remark \ref{locfin}
$$
\cn_S\ca\eqdef \tau^*(\cn_S)\otimes_{\co_C}\ca =\Ker (\ca\rightarrow \ca_{S_{red}}),
\mbox{\quad} \cn_S\ca_0 \eqdef \tau^*(\cn_S)\otimes_{\co_C}\ca_0=\Ker (\ca_0\rightarrow \ca_{S_{red},0})
$$
and
$$
\ca_{S_{red}}^*=\ca^*/(1+\cn_S\ca), \mbox{\quad} \ca_{S_{red},0}^*=\ca_0^*/(1+\cn_S\ca_0).
$$
Since $\ca /\ca_0$ is flat over $S$, we obtain, by comparing the exact sequences
$$
0\rightarrow \cn_S\ca_0\rightarrow \ca_0\rightarrow \ca_{S_{red},0}\rightarrow 0, \mbox{\quad}
0\rightarrow \cn_S\ca\rightarrow \ca\rightarrow \ca_{S_{red}}\rightarrow 0,
$$
that $\cn_S\ca\cap \ca_0=\cn_S\ca_0$ and therefore $1+\cn_S\ca_0=(1+\cn_S\ca)\cap\ca_0^*$.

Note that, by definition, the order map on $\ca^*$ coincide with the order map on $\ca^*_{S_{red}}$. Summarized we get
\begin{prop}
\label{orderfn}
\begin{enumerate}
\item \label{item1} If $\xo_{\infty}$ is a ribbon over a Noetherian scheme $S$ satisfying the assumptions in the beginning of this section, then
\begin{enumerate}
\item
\label{1))}
The order map
$$
\ord :\ca^*\lto \dz_C
$$
is a morphism of sheaves of groups.
\item
\label{2))}
There exist neigbourhoods $U_{P_s}\subset C$ such that for each $U_{P_s}$ the map $\ord |_{U_{P_s}}$ is a surjective morphism.
\item
\label{3))}
We have the equality of sheaves
$$
\Ker (\ord )=\ca_0^*\cdot (1+\cn_S\ca )\simeq \ca_0^*\coprod_{1+\cn_S\ca_0}(1+\cn_S\ca ),
$$
where on the right hand side we consider the amalgamated sum.
\end{enumerate}
\item \label{item2} If $\xo_{\infty}$ is a ribbon over a Noetherian scheme $S$ satisfying the assumptions from remark \ref{osobu}, then the statement (\ref{2))}) of this propostion  holds for $U_{P_s}=C$.
\item \label{item3} If $\xo_{\infty}$ is a ribbon obtained by the base change from a ribbon over a field $k$ of characteristic zero that satisfies the assumptions from remark \ref{osobu}, then
\begin{enumerate}
\item
\label{1)))} Using exponential and logarithmic maps we can write the equality of sheaves
$$
\Ker (\ord )=\ca_{0}^*\coprod_{\cn_S\boxtimes_{\co_S}\ca_0}\cn_S\boxtimes_{\co_S}\ca ,
$$
where on the right hand side we consider the amalgamated sum.
\item
\label{2)))}
We have the following exact sequences of sheaves
$$
1\lto \Ker (\ord )\lto \ca^*\lto \dz_C\lto 0,
$$
$$
0\lto \cn_S\boxtimes_{\co_S} \ca /\ca_0\lto \ca^*/\ca_{0}^*\lto \dz_{C}\lto 0,
$$
$$
1\lto \ca_0^*\lto \Ker (\ord )\lto \cn_S\boxtimes_{\co_S} \ca /\ca_0\lto 0.
$$
\end{enumerate}
\end{enumerate}
\end{prop}

\Proof. The proof of statements~\ref{item1},~\ref{item2} was given above. The proof of statement~\ref{item3}  follows when we use the power series for $\log (1+z)$ to identify: $1+\cn_S\ca\simeq \cn_S\ca$, $1+\cn_S\ca_0\simeq \cn_S\ca_0$ and $1+\cn_S\ca /((1+\cn_S\ca )\cap \ca_0^*)= 1+\cn_S\ca /1+\cn_S\ca_0$ with $\cn_S\ca /\cn_S\ca_0=
\cn_S \boxtimes_{\co_S} \ca /\ca_0$.

\subsection{Vanishing theorems}

In this subsection we will prove some facts which we will use later and which
may be of independent interest.

\begin{theo} \label{zar}
Let $\pi: X \to S$ be a proper morphism between schemes such that the fibres of this morphism are irreducible schemes. Then in the Zariski topology we have that the sheaf
$$
R^1\pi_* \dz = 0 \mbox{.}
$$
\end{theo}

\Proof.
We suppose that the sheaf
$R^1\pi_* \dz \neq 0 $. Then there is a point $s \in S$
such that the stalk $(R^1\pi_* \dz)_s \neq 0 $. By definition,
$$(R^1\pi_* \dz)_s =  \limind_{ U } H^1(\pi^{-1}U, \dz) =
\limind_{ U } \limind_{\{ U_{\alpha}\}_{\alpha \in I}}
\check{H}^1( \{ U_{\alpha}\}_{\alpha \in I} ,\dz)
 \mbox{,}
 $$
where $U$ runs over all open neighbourhood of the point $s$, i.e,
$S \supset U \supset s$, and  $\{ U_{\alpha}\}_{\alpha \in I}$
runs over all open covers of $\pi^{-1}(U)$, i.e.,
$\bigcup\limits_{\alpha \in I} U_{\alpha}= \pi^{-1}(U)$. Therefore there is a fixed open $U$, a fixed cover $\{ U_{\alpha}\}_{\alpha \in I}$ and
a fixed element $c \in \check{Z}^1(\{ U_{\alpha}\}_{\alpha \in I} ,\dz)$
such that the image of the element $c$ in the group $(R^1\pi_* \dz)_s$
is not equal to zero.

We define a set
$$I_0 = \{\alpha\in I \mid U_{\alpha}\cap \pi^{-1}(s) \neq \emptyset \} \mbox{.}$$
We have that $\bigcup\limits_{\alpha \in I_0} U_{\alpha} \supset \pi^{-1}(s)$.
We define a closed subset $F = X \setminus \bigcup\limits_{\alpha \in I_0} U_{\alpha}$. Then we have  $F \cap \pi^{-1}(s)=\emptyset$.
Since $\pi$ is a proper morphism, $\pi(F)$ is a closed subset in $S$,
and $s \notin \pi(F)$.

We obtain that a set $V = U \cap (S \setminus  \pi(F)) $
is an open neighbourhood of the point $s$. We have $s \in V \subset U$.
Now we have that a set
$\{ V_{\alpha} \eqdef U_{\alpha} \cap \pi^{-1}(V) \}_{\alpha \in I_0}$
is an open cover of the set $\pi^{-1}(V)$. Indeed, let a point
$x \in \pi^{-1}(V) $, and we suppose that $x \notin U_{\alpha}$  for any
$\alpha \in I_0$.
Then $x \in X \setminus  \bigcup\limits_{\alpha \in I_0} U_{\alpha} =F$.
Therefore $\pi(x) \in \pi(F) \cap V = \pi(F) \cap (U \cap (S \setminus  \pi(F))) = \emptyset$, a contradiction.

But  for any subset $J \in I_0$ we have now that
$
\bigcap\limits_{\alpha \in J} V_{\alpha} \neq \emptyset \mbox{,}
$
because
$
\bigcap\limits_{\alpha \in J} V_{\alpha} \cap \pi^{-1}(s) \neq \emptyset \mbox{,}
$
since $\pi^{-1}(s)$ is irreducible. Therefore
$\check{H}^1( \{ V_{\alpha}\}_{\alpha \in I_0} ,\dz)=0$.
The cover
$\{ V_{\alpha}  \}_{\alpha \in I_0}$ is a refinement of the cover
$\{ U_{\alpha} \cap \pi^{-1}(V)  \}_{\alpha \in I}$. Therefore the image of the element
$c \in \check{Z}^1(\{ U_{\alpha}\}_{\alpha \in I} ,\dz)$
in the group  $\check{H}^1( \{ V_{\alpha}\}_{\alpha \in I_0} ,\dz)$, and consequently in the group $(R^1\pi_* \dz)_s$
is  equal to zero. We obtained a contradiction.
\begin{flushright}
$\square$
\end{flushright}

Now we investigate the similar question in the \'{e}tale topology.

\begin{theo}
\label{etal}
Let $C$ be a normal variety over an algebraically closed field $k$ of characteristic zero, $S$ be a $k$-scheme, $X:=C\times_kS$,  $\pi :X\rightarrow S$ be the projection morphism.

Then $R^1\pi_*(\dz_{X_{et}})=0$.
\end{theo}

\Proof. Our arguments will be similar to the arguments in the proof of theorem 2.5 in \cite{We}.

It is enough to proof that for any geometric point $\bar{s}$
of $S$ $R^1\pi_*(\dz_{X_{et}})_{\bar{s}}=0$. By \cite[ch. III,
th.1.15]{Mi} we have $R^1\pi_*(\dz_{X_{et}})_{\bar{s}}\simeq
H^1_{et}(X\times_S\Spec \co_{S,s}^{sh},\dz )$. So, we can assume $S$
is a spectrum of a strict hensel ring $R$. Since every hensel ring
is the union of henselizations of its finitely generated subrings,
we can apply \cite[ch. III, lemma 1.16]{Mi}
and assume that $R$ is a strict henselization of a finitely
generated ring. By \cite[IV.18.7.3]{EGA4}, $R$ is a pseudo-geometric ring (or a universally japanese ring).

 We are going to use induction on dimension of $R$, where $R$ is a strict hensel pseudo-geometric ring.
If $\dim R=0$, then we can assume $R$ is a field,
because by \cite[th. 7.6, corol. 7.6.1]{We} $H^1_{et}(X,\dz )=H^1_{et}(X_{red},\dz )$.
Then $H^1_{et}(X,\dz )=0$  by prop. 7.4, th.7.6 in \cite{We}, because $C\times_kR$ is a normal scheme
by \cite[Ch.V, \S 1, prop. 19]{Bu} and \cite[vol.I, ch.III, \S 15, th.40, cor.1]{ZS}.

Now let $\dim R>0$. Let's consider two cases.

{\it Case 1.} $R$ is a domain.
Let $\tilde{R}$ be the normalization of $R$.
Since $R$ is a pseudo-geometric and strict hensel ring, the ring $\tilde{R}$ is also a strict hensel domain.

The scheme $\tilde{X}:=X\otimes_R\tilde{R}=C\times_k\tilde{R}$ is normal by \cite[Ch.V, \S 1, prop. 19, cor. 1]{Bu} and
\cite[vol.I, ch.III, \S 15, th.40, cor.1]{ZS}.
Let $I=Ann_R(\tilde{R}/R)$ be the conductor ideal in $R$. Then we have an isomorphism $\phi :I\rightarrow \Hom_R(\tilde{R}, R)$ by the following rule:
$\phi (i)(r)=ir$, $i\in I$, $r\in \tilde{R}$. Since $R$ is pseudo-geometric, $\tilde{R}$ is a finite $R$-module. Therefore,
we have $I\neq 0$ and $R/I$ is a strict hensel ring with $\dim R/I<\dim R$. By \cite[IV.7.7.2]{EGA4}, $R/I$ is also a pseudo-geometric ring.
Denote by $Y$ the subscheme $X\times_R(R/I)\subset X$, and by $\tilde{Y}$ the subscheme $X\times_R (\tilde{R}/I)$.
Now we are in the situation of \S 7, prop.7.8 of \cite{We}. So, we have the following long exact sequence
\begin{equation}
\label{qwerty}
0\rightarrow H^0(X,\dz )\rightarrow H^0(\tilde{X},\dz )\times H^0(Y, \dz )
\rightarrow H^0(\tilde{Y}, \dz )\rightarrow H^1_{et}(X, \dz )
\rightarrow H^1_{et}(\tilde{X},\dz )\times H^1_{et}(Y, \dz ).
\end{equation}
By induction on dimension of the ring $R$ we have $H^1_{et}(Y, \dz )=0$. $H^1_{et}(\tilde{X},\dz )=0$
by prop. 7.4, th.7.6 in \cite{We}, because $\tilde{X}$ is a normal scheme. Since $\tilde{X}$, $\tilde{Y}$ are connected schemes,
the map $H^0(\tilde{X},\dz )\rightarrow H^0(\tilde{Y}, \dz )$ is surjective. Therefore, $H^1_{et}(X, \dz )=0$.

{\it Case 2.} In general case, we can assume that $R=R_{red}$ by \cite[corol. 7.6.1]{We}.
Let $(0)=\wp_1\cap\ldots\cap\wp_n$ be a primary decomposition in $R$. Set $\tilde{R}=R/\wp_1\times\ldots\times R/\wp_n$.
Set $I=\bigoplus\limits_i\bigcap\limits_{j\neq i}\wp_j$.
We have
$$
I=\bigoplus_i\bigcap_{j\neq i}\wp_j=\bigoplus_i\Hom\nolimits_R(R/\wp_i,R)
$$
is a conductor ideal in $R$, and it contains a nonzerodivisor.
Now in the notations of case 1 we have
$\tilde{X}=\coprod\limits_{i=1}^n(X\times_RR/\wp_i)$,
$\tilde{Y}=\coprod\limits_{i=1}^n(X\times_RR/(\wp_i+\cap_{j\neq i}\wp_j))$.
So, in the sequence (\ref{qwerty}) we have $H^1_{et}(Y, \dz )=0$ by induction on dimension of the ring $R$,
since $\dim R/I<\dim R$; $H^1_{et}(\tilde{X},\dz )=\prod\limits_{i=1}^n H^1_{et}(X\times_RR/\wp_i, \dz )=0$ by case 1;
and $H^0(\tilde{X},\dz )\rightarrow H^0(\tilde{Y}, \dz )$ is a surjective map. So, again $H^1_{et}(X, \dz )=0$.
\begin{flushright}
$\square$
\end{flushright}

Now we can prove the same result in the flat topology, where under the flat topology we understand the fppf or fpqc topology on $X$.

\begin{theo}
\label{flat1}
Let $C$ be a normal variety over an algebraically closed field $k$ of characteristic zero, $S$ be a $k$-scheme,
$X:=C\times_kS$,  $\pi :X\rightarrow S$ be the projection morphism.

Then $R^1\pi_*(\dz_{X_{fl}})=0$.
\end{theo}

\Proof. By definition of the sheaf $R^1\pi_*(\dz_{X_{fl}})$, we need to prove that for any element $x\in H^1_{fl}(C\times_k U, \dz )$, where $U$ is flat over $S$, there exists a cover $(U_i)$ of $U$ in the flat topology such that $res_{U,U_i}(x)=0$ for any $i$.

For any scheme $Z$ set $Y=Z[t,t^{-1}]=Z\times_{\sdz}\dz [t,t^{-1}]$, and let $p:Y\rightarrow Z$ be the structure map. Similarly, let $p^+$ and $p^-$ denote the structure maps from $Y^+=Z[t]=Z \times_{\sdz}\dz [t]$ and $Y^-=Z[t^{-1}]=Z \times_{\sdz}\dz [t^{-1}]$ to $Z$.

Recall that for a covariant functor $F$ from the category of
commutative rings (or for a contravariant functor from the category
of schemes) to some Abelian category the following functors are
defined (see \cite{We}, \S 1 or \cite{B}, ch.XII):
$$
NF(R)=N_tF(R)=\Ker [F(t=1):F(R[t])\rightarrow F(R)]\simeq \Coker [F(i_+):F(R)\rightarrow F(R[t])];
$$
$$
LF(R)=\Coker [F(R[t])\oplus F(R[t^{-1}])\overset{add}{\lto} F(R[t,t^{-1}])].
$$
Clearly, $F(R[t])\simeq F(R)\oplus NF(R)$.

Consider now the following sequence of sheaves on $Z$ (in the flat topology) from \cite{We}, proof of prop. 7.2:
\begin{equation}
\label{units}
0\lto \G_{m,Z}\lto p_*^+(\G_{m,Y^+})\times p_*^-(\G_{m,Y^-})\lto p_*(\G_{m,Y})\lto \dz \lto 0.
\end{equation}
By (1.1) in \cite{We} this sequence is exact and there is a splitting $\dz \rightarrow p_*(\G_{m,Y})$ given by multiplication by $t$ ($n\mapsto t^n$) and the splitting $p_*(\G_{m,Y})\rightarrow \G_{m,\sdz }$ given by evaluation at $t=1$. So, we have
\begin{equation}
\label{isom}
p_*^+(\G_{m,Y^+})\simeq \G_{m,Z}\times \cn_t\G_{m,Z}, \mbox{\quad } p_*^-(\G_{m,Y^-})\simeq \G_{m,Z}\times \cn_{t^{-1}}\G_{m,Z}
\end{equation}
and
$$
p_*(\G_{m,Y})\simeq \G_{m,Z}\times \cn_t\G_{m,Z}\times \cn_{t^{-1}}\G_{m,Z}\times\dz ,
$$
where $\cn_t\G_{m,Z}$ is the sheaf, associated to the presheaf $U\mapsto \Coker(\G_{m,Z}(U)\rightarrow p^+_*(\G_{m,Y^+})(U))$, in the flat topology.

Now, comparing  the Leray spectral sequences for $p^+$, $p^-$ and
$p$ (in the flat topology), we obtain the following exact diagram
(compare the diagram in the proof of th.7.6 in \cite{We}):
$$
\begin{array}{ccccccccc}
&&&0&&0&&&\\
&&&\downarrow &&\downarrow &&&\\
0\rightarrow & Pic(Z)&\rightarrow &H^1_{fl}(Z, p^+_*(\G_{m,Y^+}))\oplus &\rightarrow &H^1_{fl}(Z, p_*(\G_{m,Y}))&\rightarrow &H^1_{fl}(Z,\dz )&\rightarrow 0\\
&&&H^1_{fl}(Z, p^-_*(\G_{m,Y^-}))&&&&&\\
&\parallel&&\downarrow &&\downarrow &&\downarrow &\\
0\rightarrow & Pic(Z)&\rightarrow &Pic(Y^+)\oplus Pic(Y^-)&\rightarrow & Pic(Y)&\rightarrow &LPic(Z)& \rightarrow 0\\
&&&\downarrow &&\downarrow &&&\\
&0&\rightarrow &H^0(Z,R^1p^+_*(\G_{m,Y^+}))\oplus &\rightarrow &H^0(Z,R^1p_*(\G_{m,Y}))&&\\
&&&H^0(Z,R^1p^-_*(\G_{m,Y^-}))&&&&&
\end{array}
$$
Here the first row is exact by (\ref{units}) and (\ref{isom}).
In the second row we use Hilbert's 90 theorem: $H^1_{fl}(T,\G_{m,T})=Pic (T)=H^1_{et}(T, \G_{m,T})$, where $T$ is a scheme. So, the second row is exact by \cite[th. 7.6]{We}. Let's show that the third row is also exact.

Assume the converse, and let $x\in H^0(Z,R^1p^+_*(\G_{m,Y^+}))$ be an element from the kernel of the map from the third row. This means that there exists a cover $(U_i)$ of $Z$ in the flat topology such that
$$res_{Z,U_i}(x)\in \Ker (H^1_{fl}(U_i\times_{\sdz} \dz [t], \G_m)\rightarrow H^1_{fl}(U_i\times_{\sdz} \dz [t, t^{-1}], \G_m))$$
for all $i$. By Hilbert's 90 theorem we have
$$res_{Z,U_i}(x)\in \Ker (H^1_{et}(U_i\times_{\sdz} \dz [t], \G_m)\rightarrow H^1_{et}(U_i\times_{\sdz} \dz [t, t^{-1}], \G_m)).$$

By $(5.1)$ in \cite{We} (see also the proof of th. 7.6 there) we have $R^1p^+_*(\G_{m,Y^+})\oplus R^1p^-_*(\G_{m,Y^-})\simeq R^1p_*(\G_{m,Y})$ in the \'etale topology on $Z$ for any scheme $Z$. So, our conditions mean that for every $i$ there exists a cover $(V_{i\alpha })$ of $U_i$ such that $0=res_{U_i,V_{i\alpha }}(res_{Z,U_i}(x))\in H^1_{et}(V_{i\alpha }\times_{\sdz} \dz [t], \G_m)$. Again by Hilbert's 90 theorem this means $res_{Z,V_{i\alpha }}(x)=0$ for the flat cover $(V_{i\alpha })$ of $Z$, i.e. $x=0$. The same arguments work for $x\in H^0(Z,R^1p^-_*(\G_{m,Y^-}))$.

A diagram chase now shows that $H^1_{fl}(Z, \dz )\rightarrow LPic(Z)\simeq H^1_{et}(Z, \dz )$ is an injective map. Note also that for an open \'etale $U\rightarrow Z$ the diagram
$$
\begin{array}{ccc}
H^1_{fl}(Z, \dz )&\rightarrow &LPic(Z)\\
\downarrow &&\downarrow \\
H^1_{fl}(U, \dz )&\rightarrow &LPic(U)
\end{array}
$$
is commutative.

Now let $x\in H^1_{fl}(C\times_k U, \dz )$, where $U$ is flat over $S$. By the arguments above, the element $x$ is embedded in $H^1_{et}(C\times_k U, \dz )$. By theorem~\ref{etal} there exists an \'etale  cover $(U_i)$ of $U$ such that $0=res_{U,U_i}(x)\in H^1_{et}(C\times_k U_i, \dz )$. Then by the arguments above this imply that for the cover $(U_i)$ in the flat topology we also have $0=res_{U,U_i}(x)\in H^1_{fl}(C\times_k U_i, \dz )$. So, $R^1\pi_*(\dz_{X_{fl}})=0$.
\begin{flushright}
$\square$
\end{flushright}

\subsection{Representability of the Picard functor of $X_{\infty}$}
\label{repxinf}
Let $\xo_{\infty}= (C, \ca)$ be a
ribbon over a field $k$.
Then we get the locally ringed space $X_{\infty}=(C, \ca_0)$.

\begin{defin} \label{defpicrib}
We denote by ${\bf {Pic}}_{X_{\infty}}$ the sheaf
   on
the big Zariski site of $\Spec k$, associated with the functor $ S
\mapsto {Pic}_{X_{\infty}}(S)$ (In other words, for any Noetherian
scheme $S$ over $k$ we consider all scheme-theoretic open affine
covers of $S$ and we take the sheaf associated with the presheaf
$S'\mapsto Pic_{X_{\infty}}(S')$ with respect to these covers. ).
\end{defin}

\begin{rem}
From definitions~\ref{defpicrib} and~\ref{picarf}  it follows that for any Noetherian scheme $S$
we have
$$
{\bf {Pic}}_{X_{\infty}}(S)= H^0(S,R^1\pi_* \ca_{S,0}^*)  \mbox{,}
$$
where $\pi : C\times_k S \to S$ is the projection morphism.
\end{rem}

\begin{rem}
\label{cites}
If the curve  $C$ is proper over $k$, then the locally ringed space $X_{\infty}=(C,\ca_0)$ is a
weakly Noetherian formal scheme in the sense of \cite{Lip}. Then for the
field $k$ of any characteristic Lipman proved  in~\cite[section 2.5]{Lip} that the fpqc
sheaf associated with the modified Picard functor of $X_{\infty}$ is
 a $k$-group scheme.

Under assumptions that  $\chara k = 0$, the field $k$ is algebraically closed,
and  $C$ is a projective irreducible curve we will give now an easy proof that the sheaf
${\bf {Pic}}_{X_{\infty}}$ is a $k$-group scheme. We will study also the structure of this $k$-group scheme. We note that from  the existence of this
$k$-group scheme it will automatically follows that the presheaf ${\bf {Pic}}_{X_{\infty}}$ is a sheaf on the big fpqc site of $\Spec k$.
\end{rem}

We will need the following lemma and corollaries of this lemma.
\begin{lemma}
\label{lamma}
Let $\xo_{\infty}= (C, \ca)$ be a
ribbon over a field $k$.
Let $S$ be an affine scheme over the field $k$, and $\cm$ be a coherent sheaf
on $S$.
We have
$$
H^h(C\times_k S, \ca_i\widehat{\boxtimes}_k \cm)\simeq \limproj_{j>i}H^h(C\times_k S, (\ca_i/\ca_j){\boxtimes}_k \cm)\simeq \limproj_{j>i}(H^h(C, \ca_i/\ca_j)\otimes_kH^0(S, \cm)),
$$
$$
H^q(C\times_k S, \ca_i\widehat{\boxtimes}_k \cm)=0
$$
for any $i\in \dz$, $h\le 1$, $q>1$.
\end{lemma}

\Proof. We have the analog of the K\"unneth formula:
\begin{equation} \label{kf}
p_*((\ca_i/\ca_{i+h}){\boxtimes}_k \cm)\simeq
(\ca_i/\ca_{i+h})\otimes_kH^0(S,\cm),
\end{equation}
where $p :C\times_k S\rightarrow C$ is the projection. Indeed, if
$U$ is an affine open set on $C$, and $\tau_S: C\times_k
S\rightarrow S$ is the projection, we have the natural isomorphisms
$$
p_*((\ca_i/\ca_{i+h}){\boxtimes}_k \cm)(U)\simeq
p^*(\ca_i/\ca_{i+h})(U\times_k S) \otimes_{\co_{U \times_k S}}
\tau_S^*\cm(U\times_kS)\simeq
$$
$$
\simeq  (\ca_i/\ca_{i+h})(U)\otimes_kH^0(S,\cm),
$$
since $(\ca_i/\ca_{i+h})$, $\cm$ are coherent sheaves of modules
on $X_{h-1}$, $S$ correspondingly (see prop.~1 in \cite{Ku}). These
isomorphisms are obviously compatible with the restriction
homomorphisms corresponding to the embedding of affine sets
$U'\subset U$ for both sheaves $p_*((\ca_i/\ca_{i+h}){\boxtimes}_k
\cm)$ and $(\ca_i/\ca_{i+h})\otimes_kH^0(S,\cm)$. Therefore, the
sheaves from formula~(\ref{kf}) are isomorphic.

Since $p$ is an affine morphism, we have then
$$H^q(C\times_k S, (\ca_i/\ca_{i+h}){\boxtimes}_k \cm )\simeq H^q(C, \ca_i/\ca_{i+h})\otimes_k H^0(S,\cm)$$
for all $q$ (see \cite[ch.III, ex.8.1]{Ha}).

For all $i,h,k$ with $h\le k$ we have surjective morphism of sheaves
$$ (\ca_i/\ca_{i+k}){\boxtimes}_k \cm \rightarrow
(\ca_i/\ca_{i+h}){\boxtimes}_k \cm \mbox{,}$$ because
$(\ca_i/\ca_{i+k}) \rightarrow (\ca_i/\ca_{i+h})$ is a surjective
morphism of sheaves on $C$.

For any affine $U$ in $C\times_k S$ the maps $$\Gamma (U,
(\ca_i/\ca_{i+k}){\boxtimes}_k \cm)\rightarrow \Gamma (U,
(\ca_i/\ca_{i+h}){\boxtimes}_k \cm)$$ are surjective, since
$(\ca_l/\ca_m) {\boxtimes}_k \cm$ are coherent sheaves of modules
for all $l<m$ on $X_{m-l-1}\times_kS$. By the same reason we have
$H^q(U, (\ca_l/\ca_m) {\boxtimes}_k \cm)=0$ for all $q>0$.

At last, since $C$ is a projective curve, the projective systems
$$\{H^q(C, \ca_i/\ca_{i+h})\otimes_k H^0(S,\cm)\}_{h\in\sdn}
\mbox{,} \quad q\ge 0$$
satisfy the ML-condition. So, by \cite[ch.
0, prop.13.3.1]{EGA3} we have
$$
H^q(C\times_k S, \ca_i\widehat{\boxtimes}_k \cm)\simeq \limproj_{j>i}H^q(C\times_k S, (\ca_i/\ca_j){\boxtimes}_k \cm)\simeq \limproj_{j>i}(H^q(C, \ca_i/\ca_j)\otimes_kH^0(S, \cm))
$$
for $q\ge 1$.

For $q=0$ it follows from the definition of the sheaf $\ca_i\widehat{\boxtimes}_k \cm$.
\begin{flushright}
$\square$
\end{flushright}

\begin{corol}
\label{corol1}
For a ribbon $\xo_{\infty}= (C, \ca)$  over a field $k$, for
 $S$ which is an affine Noetherian scheme over the field $k$, and for $\cm$
which is  a coherent sheaf
on $S$  we have
$$H^1(C\times_k S, (\ca /\ca_0)\widehat{\boxtimes}_k\cm)\simeq H^1(C, \ca /\ca_0) \otimes_kH^0(S, \cm)\mbox{\quad }
$$
$$
H^q(C\times_k S, \ca\widehat{\boxtimes}_k\cm)=H^q(C\times_k S, (\ca /\ca_0){\boxtimes}_k\cm)=0$$
for $q\ge 2$.
\end{corol}

\Proof. The proof is clear, since cohomology commute with $\limind$
on Noetherian schemes.
\begin{flushright}
$\square$
\end{flushright}

\begin{corol}
\label{corola1}
For a ribbon $\xo_{\infty}= (C, \ca)$  over an algebraically closed field $k$
with a projective irreducible curve $C$, for
 $S$ which is an affine  scheme over the field $k$
there is an embedding of $k$-algebras
$$H^0(C \times_k S, \co_{C \times_kS}) \lto
 H^0(C\times_kS, \ca_{S,0})
$$
 which splits the natural map of $k$-algebras
$$
H^0(C\times_kS, \ca_{S,0}) \lto H^0(C \times_k S, \co_{C \times_kS}) \mbox{.}
$$
\end{corol}
\Proof.
 By formula~(\ref{kf}),  and since $H^0(C,\co_C)\simeq k$, we have
$$H^0(C\times_kS, \co_{C \times_kS})\simeq H^0(C,\co_C)\otimes_kH^0(S,\co_S)\simeq H^0(S,\co_S) \mbox{.}$$
 Now since we have
embeddings
$$H^0(S,\co_S)\hookrightarrow H^0(C,
\ca_0/\ca_j)\otimes_kH^0(S,\co_S)= H^0(C\times_k S, \ca_{S,j}/\ca_{S,0} )$$
for all $j>0$, we obtain the
embedding
$$H^0(S,\co_S)\hookrightarrow H^0(C\times_kS, \ca_{S,0})$$
for any scheme $S$ by lemma \ref{lamma}.
\begin{flushright}
$\square$
\end{flushright}

\begin{prop} \label{reprx}
Let  $\xo_{\infty}= (C, \ca)$ be a ribbon  over an algebraically closed field $k$, $\chara k = 0$, and $C$ be
 a projective irreducible curve. We have the following properties.
\begin{enumerate}
\item The sheaf ${\bf {Pic}}_{X_{\infty}}$  is a $k$-group scheme.
\item The following sequence of $k$-group schemes is exact:
\begin{equation} \label{exseq}
 0 \arrow{e} \dv \arrow{e}  {\bf {Pic}}_{X_{\infty}}
\arrow{e,t}{\phi} {\bf {Pic}}_C  \arrow{e} 0  \mbox{,}
\end{equation}
where $\dv$ is an affine $k$-group scheme, and ${\bf {Pic}}_C$
is the Picard variety of the curve $C$, whose connected component is the generalized Jacobian of the curve $C$. There is a splitting of the
map~$\phi$ from sequence~(\ref{exseq})
over any affine subscheme $U$ of the scheme ${\bf {Pic}}_C$.
\end{enumerate}
\end{prop}

\Proof.
Since we supposed that $\chara k =0$, we can use the
series for $exp(z)$ and $log(1+z)$.

For any affine Noetherian scheme $S$ we have exact sequences of
sheaves on $C \times_k S$:
$$
0\rightarrow \ca_{S,1}\overset{exp}{\rightarrow}
\ca_{S,0}^*\rightarrow \co_{C \times_kS}^*\rightarrow 1 \mbox{.}
$$
Therefore, using lemma~\ref{lamma} and corollary~\ref{corola1} of this lemma, we
obtain the following exact sequence
\begin{equation} \label{pbun}
0\rightarrow H^1(C, \ca_1)\widehat{\otimes}_kH^0(S,\co_S)\rightarrow Pic_{X_{\infty}}(S)\rightarrow Pic_C(S)\rightarrow 0 \mbox{,}
\end{equation}
where $H^1(C, \ca_1)\widehat{\otimes}_kH^0(S,\co_S)\eqdef
\limproj\limits_{j>1}(H^1(C, \ca_1/\ca_j)\otimes_k H^0(S, \co_S))$,
and
$$S \mapsto Pic_C(S) = H^1(C \times_kS, \co_{C \times_k S}^*)$$
is the Picard functor  of the curve $C$.

We define
$$
H\eqdef \Hom\nolimits_{k,cont}(H^1(C,\ca_1),k)=\limind_{j>1}(H^1(C,
\ca_1/\ca_j))^* \mbox{,}
$$
and $\dv \eqdef \Spec (Sym_k(H))$ is an affine $k$-group scheme,
where $Sym_k(H) \eqdef \bigoplus\limits_{i =0}^{\infty} S^i(H)$, and
the group law is given by $v \mapsto v \otimes 1 + 1 \otimes v$, $v
\in H$.

We have for any affine Noetherian scheme $S$ over $k$
\begin{multline}
\Hom\nolimits_{sch}(S, \dv )=\Hom\nolimits_{k-alg}(Sym_k(H),H^0(S,\co_S))=\Hom\nolimits_k(H, H^0(S,\co_S))= \\
\limproj_{j>1}(H^1(C, \ca_1/\ca_j)\otimes_kH^0(S, \co_S))=H^1(C, \ca_1)\widehat{\otimes}_kH^0(S,\co_S).
\end{multline}
Thus, we have from exact sequence~(\ref{pbun}) the following exact sequence of groups, which is functorial with
$S$:
\begin{equation} \label{lastseq}
0\arrow{e} \Hom\nolimits_{sch}(S, \dv ) \arrow{e}
Pic_{X_{\infty}}(S)
\arrow{e,t}{\phi} Pic_C(S)\rightarrow 0,
\end{equation}
The Zariski sheaves associated with the first and the last presheaves (or
functors) in sequence~(\ref{lastseq}) are $k$-group schemes,
because for the first sheaf it is true by construction, and the last sheaf
is the Picard scheme ${\bf {Pic}}_C$ of the curve $C$ (see~\cite{Gr}).

Representability of ${\bf {Pic}}_C$ means that there is a universal
object $\lambda$ (a Poincar\'e bundle on $C\times_k{\bf {Pic}}_C$)
corresponding to the identity map of ${\bf {Pic}}_C$ under $\Hom
({\bf {Pic}}_C, {\bf {Pic}}_C)\simeq {\bf {Pic}}_C({\bf {Pic}}_C)$
(we do not distinguish here in the notation for a representable
functor and the scheme that represents it). By the construction of
the associated sheaf, $\lambda$ is given by an open (affine)
covering $\{U_{\alpha}\}$ of ${\bf {Pic}}_C$, and line bundles
$\cl_{\alpha}$ on $C\times_kU_{\alpha}$ ($\cl_{\alpha}$,
$\cl_{\beta}$ are isomorphic on $C\times_k(U_{\alpha}\cap
U_{\beta})$ up to twist with a line bundle on $U_{\alpha}\cap
U_{\beta}$). The functorial map $\Hom (\cdot ,U_{\alpha})\rightarrow
{\bf {Pic}}_C(\cdot )$ is given by the embedding $U_{\alpha} \subset
{\bf {Pic}}_C$.

Since by (\ref{lastseq}) the bundles $\cl_{\alpha}$ on $C\times_kU_{\alpha}$ can be lifted to line bundles $\widetilde{\cl_{\alpha}}$ on $X_{\infty ,U_{\alpha}}$, the bundles $\widetilde{\cl_{\alpha}}$ give rise to morphisms of functors $s_{\alpha}: \Hom (\cdot ,U_{\alpha})\rightarrow {\bf {Pic}}_{X_{\infty}}(\cdot )$, which, composed with $\phi : {\bf {Pic}}_{X_{\infty}} \rightarrow {\bf {Pic}}_C$, gives the embedding $U_{\alpha} \subset {\bf {Pic}}_C$.

Thus, as a functor we can split $P_{\alpha}\eqdef {\bf {Pic}}_{X_{\infty}}\times_{{\bf {Pic}}_C}U_{\alpha} \subset {\bf {Pic}}_{X_{\infty}}$ into $\dv\times U_{\alpha}$ by the action of $\dv$ on $P_{\alpha}$ induced from the group structure ($P_{\alpha}\rightarrow \dv\times U_{\alpha}$ by $x\mapsto (x-s_{\alpha}\phi (x), \phi (x))$, $\dv\times U_{\alpha}\rightarrow P_{\alpha}$ by $(x', x_0)\mapsto x'+s_{\alpha}(x_0)$).

Thus, ${\bf {Pic}}_{X_{\infty}}$ has a cover $\{P_{\alpha}\}$ by
open, representable functors, hence is representable.

We consider any affine subscheme $U$ of the scheme ${\bf {Pic}}_C$.
Then restricting a Poincar\'e bundle from $C\times_k{\bf {Pic}}_C$
to $U \times_k{\bf {Pic}}_C$ and using the same arguments as above
it is easy to see that there is a splitting of the map $\phi$ from
sequence~(\ref{exseq}) over the subscheme $U$ of the scheme ${\bf
{Pic}}_C$. The proposition is proved.

\begin{flushright}
$\square$
\end{flushright}

\begin{rem}
Suppose we have that the following property is satisfied:
${\tau_{S}}_*(\ca_{S,0})=\co_S$, where $\tau_S :X_{\infty
,S}\rightarrow S$ is a structure morphism (between the locally
ringed spaces). This condition is satisfied, for example, if the
ribbon $\xo_{\infty}$ comes from a surface and an ample Cartier
divisor, because in this case $H^0(C, \ca_1/\ca_i)=0$ for any $i \ge
1$ (see remark \ref{r2} and theorem~\ref{pc}), where from, using the
exact sequence
$$
0\rightarrow \ca_{S,1} \rightarrow \ca_{S,0} \rightarrow \co_{C\times_kS}\rightarrow 0,
$$
lemma \ref{lamma} and corollary \ref{corola1}, we obtain
${\tau_{S}}_*(\ca_{S,0})=\co_S$.

Now, since the curve $C$ has a $k$-rational point, there exists a
section $\sigma_S:S\rightarrow X_{\infty ,S}$ of the morphism
$\tau_S$. So, using standard arguments (see, for example, \cite{Gr})
we obtain for this case the following description of the functor
${\bf {Pic}}_{X_{\infty}}$:
$$
{\bf {Pic}}_{X_{\infty}}(S)\simeq Pic_{X_{\infty}}(S)/Pic (S)\simeq
\{ \cl\in Pic_{X_{\infty}}(S)| \mbox{\quad} \sigma_S^*(\cl )\simeq \co_S\}.
$$
\end{rem}

\subsection{Representability of the Picard functor of $\xo_{\infty}$}
\label{rpf}

In this subsection we will show that under certain conditions on a
ribbon $\xo_{\infty}$ the Zariski sheaf, associated with the Picard functor
of a ribbon $\xo_{\infty}$, is
representable by a formal scheme.

Let $C$ be an irreducible projective curve over an algebraically
closed field $k$ of characteristic zero, and $\xo_{\infty}$ be a
ribbon over $k$ with underlying topological space $C$ and either
with a smooth point in the sense of definition~9 from \cite{Ku}, or
satisfying the  condition $(**)$ from definition~\ref{uslovie}. We
consider the ribbon $\xo_{\infty ,S}$ for some base change
$S\rightarrow \Spec k$.

Let $\cf \in Pic(\xo_{\infty ,S})$. We define a sheaf of generating
sections $\cb (\cf )$ (which is a sheaf of sets)  by the rule
$$
\cb (\cf )(U)=\{\mbox{ sections $\lambda\in \cf (U)$ with $\cf
|_U=\ca |_U\cdot \lambda$}\} \mbox{,}
$$
where $U$ is open in $C \times_k S$.  We have $\cb (\cf
)(U)=\varnothing$ or (after choice of one generator) $\cb (\cf )|_U
\simeq \ca^*|_U$. Thus, $\cb (\cf )$ is a torsor over the sheaf of
groups $\ca^*$.

We recall that for any topological space $Y$ we denoted by $\dz_Y$
the sheaf of locally constant functions on $Y$ with values in $\dz$.
\begin{defin}
\label{def12'} Let $\cb$ be a category of affine Noetherian
$k$-schemes. We define the  contravariant  functor
$\widetilde{Pic}'_{\sxo_{\infty}}$ from $\cb$ to the category of
Abelian groups:
$$
\widetilde{Pic}'_{\sxo_{\infty}}(S)\eqdef \{\mbox{the group of
isomorphism classes of pairs $(\cf ,d)$}\} ,
$$
where $\cf\in Pic(\xo_{\infty ,S})$, and  $d:\cb (\cf )\rightarrow
\dz_{C\times_k S}$ is a morphism of sheaves of sets such that
$$
d(a\lambda )(x)=\ord (a)(x)+d(\lambda )(x)
$$
for any $a\in \ca_S^*(U)$, $\lambda\in \cb(\cf) (U)$, $U\subset
C\times_k S$ is open, and $x\in U$. Two pairs $(\cf ,d)$ and $(\cf
',d')$ are isomorphic, if there is an isomorphism of sheaves of
$\ca_S$-modules $\cf$ and $\cf '$ compatible with $d, d'$. Besides
$$
(\cf_1, d_1) \otimes (\cf_2, d_2) = (\cf_1 \otimes_{\ca_S} \cf_2, d)
\mbox{,}
$$
where $d(\lambda_1 \otimes \lambda_2) = d_1(\lambda_1) +
d_2(\lambda_2)$.
\end{defin}

\begin{ex}
Any locally free sheaf of rank $1$ on $\xo_{\infty}$ has a filtration (see example~11 in \cite{Ku}), i.e.,
this sheaf is a torsion free sheaf on $\xo_{\infty}$ in the sense of definition~\ref{tfsh}. Its base change
give some sheaf $\cf \in Pic(\xo_{\infty ,S})$, which also has a filtration.
 We define a morphism of order $d:\cb (\cf )\rightarrow W_{C\times_k S}(\dz )$  (where $W_{C\times_k S}$ is
 the sheaf of all
 functions on  $C\times_k S$ with values in $\dz$) by the rule
$$
d(\lambda )(x)= \max\{j |\mbox{\quad } \lambda |_{U_s}\in (\cf_s)_{x,j}\},
$$
where $\lambda\in \cb(\cf) (U)$ for an open $U\subset C\times_k S$, $x\in U$, $s=\tau (x)$, $U_s$ is an open set in $C_s$, which is obtained from $U$ by the base change $s\rightarrow S$, $\cf_s$ is a sheaf after the base change.
 One can easily check that $d$ is a morphism of sheaves of sets. Moreover, by section~\ref{funorder} it factors through the subsheaf
 $\dz_{C\times_k S} \subset W_{C\times_k S}$.
Besides, we have
$$d(a\lambda )(x)=\ord (a)(x)+d(\lambda )(x)$$
for any $a\in \ca_S^*(U)$, because $\ord$ is a morphism of sheaves of groups by proposition \ref{orderfn}.
 Thus, $(\cf ,d) \in \widetilde{Pic}'(\xo_{\infty ,S})$.
\end{ex}

We define the following sheaf of groups on $C \times_k S$:
\begin{equation}
\label{Ims} \Im_S \eqdef \Ker (\ord :\ca^*_S\rightarrow
\dz_{C\times_k S})=
\ca_{S,0}^*\coprod_{\cn_S\boxtimes_{k}\ca_0}\cn_S\boxtimes_{k}\ca
\mbox{,}
\end{equation}
where the last equality follows from proposition \ref{orderfn}.

\begin{defin}
\label{def12} Let $\cb$ be a category of affine Noetherian
$k$-schemes. We define the  contravariant  functor
$\widetilde{Pic}_{\sxo_{\infty}}$ from $\cb$ to the category of
Abelian groups:
$$
\widetilde{Pic}_{\sxo_{\infty}}(S)\eqdef H^1(C\times_k S, \Im_S)
 \mbox{,}
$$
where the restriction maps of this functor are compositions of the
natural maps
$$H^1(C\times_k S, \Im_S) \rightarrow H^1(C\times_k S,
(id\times f)_*(\Im_{S'}))\rightarrow H^1(C\times_k S', \Im_{S'}),
$$
and the second map is the embedding from the Cartan-Leray spectral
sequence for a morphism $f:S'\rightarrow S$.
\end{defin}

We have always an evident morphism of functors:
$$
\widetilde{Pic}_{\sxo_{\infty}} \lto
\widetilde{Pic}'_{\sxo_{\infty}}
$$
such that for any $S \in \Ob(\cb)$ we have an embedding of Abelian
groops:
$$
\widetilde{Pic}_{\sxo_{\infty}}(S) \hookrightarrow
\widetilde{Pic}'_{\sxo_{\infty}}(S) \mbox{.}
$$

\begin{prop}
\label{compare} Let a ribbon $\xo_{\infty}$ satisfies the
condition~$(**)$ from definition~\ref{uslovie}. Then  we have the
natural isomorphism of functors:
$$
\widetilde{Pic}_{\sxo_{\infty}} \simeq
\widetilde{Pic}'_{\sxo_{\infty}} \mbox{.}
$$
\end{prop}

\Proof. Let $S \in \Ob(\cb)$. The sheaf  of automorphisms of a pair
$(\cf ,d) \in \widetilde{Pic}'_{\sxo_{\infty}}(S)$ (i.e., the sheaf
of automorphisms of $\ca$-module $\cf$, which preserve the function
$d$)  is equal to the sheaf $\Im_S$. Besides, by remark~\ref{osobu},
the pair $(\cf ,d)$ is isomorphic to the pair $(\ca_S, \ord)$
locally on $C \times_k S$. Therefore, by standard arguments with
twisted forms (see, e.g. \cite{Mi}, ch.III, \S 4) we obtain the
statement of the proposition.
\begin{flushright}
$\square$
\end{flushright}

\begin{defin} \label{dddd}
We denote by ${\bf \widetilde{Pic}}_{\sxo_{\infty}}$ the sheaf on
the big Zariski site of $\Spec k$,
 associated with the presheaf $
S \mapsto \widetilde{Pic}_{\sxo_{\infty}}(S)$ (compare with
definition~\ref{defpicrib}).

Analogously, we denote by ${\bf {Pic}}_{\sxo_{\infty}}$ the sheaf
on
the big Zariski site of $\Spec k$,  associated with the presheaf $
S \mapsto {Pic}_{\sxo_{\infty}}(S)$.
\end{defin}

\begin{rem}
From definitions~\ref{dddd},~\ref{def12} and~\ref{picarf}   it follows that for any Noetherian scheme $S$
we have
$$
{\bf \widetilde{Pic}}_{\sxo_{\infty}}(S)= H^0(S,R^1\pi_* \Im_S) \qquad
\mbox{,} \qquad
{\bf \widetilde{Pic}}_{\sxo_{\infty}}(S)=H^0(S, R^1\pi_* \ca_{S,}^*) \mbox{,}
$$
where $\pi : C\times_k S \to S$ is the projection morphism.
\end{rem}

In view of theorem \ref{th1}, propositions~\ref{freeness}
and~\ref{compare}, it is important to obtain that the sheaf ${\bf
\widetilde{Pic}}_{\sxo_{\infty}}$ is a $k$-group scheme.  Our first
aim is to prove this under some conditions, and then we will compare
the sheaf ${\bf \widetilde{Pic}}_{\sxo_{\infty}}$ with the sheaf
${\bf Pic}_{\sxo_{\infty}}$.

\begin{theo}
\label{represent} Let $C$ be an irreducible projective curve over an
algebraically closed field $k$ of characteristic zero, and
$\xo_{\infty}$ be a ribbon with underlying topological space $C$,
which satisfies conditions from the beginning of section~\ref{rpf}.
Assume that
\begin{equation}
\label{assumpt} \Coker (H^0(C, \ca )\lto H^0(C, \ca /\ca_0)) =0
\mbox{.}
\end{equation}
Then the  sheaf  ${\bf \widetilde{Pic}}_{\sxo_{\infty}}$ is a formal
group scheme, which is isomorphic (non-canonically) to the product
${\bf Pic}_{X_{\infty}}\times_k{\bf \widehat{Br}_{\sxo_{\infty}}}$,
where ${\bf Pic}_{X_{\infty}}$ is the Picard scheme of $X_{\infty}$
(see proposition~\ref{reprx}), and ${\bf \widehat{Br}_{\sxo_{\infty}}}$
is the formal Brauer group of $\xo_{\infty}$ (see proposition~\ref{prop3}) .
\end{theo}
\begin{rem}
Compare assumption formula~(\ref{assumpt}) with theorem~\ref{pc}.
\end{rem}

\Proof. It is enough to prove that for any affine Noetherian scheme $S$  over $k$ the following sequence is split exact (see corollary \ref{corol1} for the last term):
\begin{equation}
\label{splitex} 0\lto H^1(C\times_kS, \ca_{S,0}^*)\lto
H^1(C\times_kS, \Im_S)\lto H^1(C, \ca /\ca_0)\otimes_kH^0(S,
\cn_S)\lto 0,
\end{equation}
and the splitting is functorial with $S$.
(We recall that the coherent sheaf $\cn_S$ is  the nilradical of the scheme $S$).
 Indeed, by~proposition~\ref{reprx} the Zariski sheaf, associated with the presheaf $S\mapsto H^1(C\times_kS, \ca_0^*)$, is a scheme ${\bf Pic}_{X_{\infty}}$.  By remark \ref{fb1} we have
$$\Hom\nolimits_{form. sch.} (S, {\bf \widehat{Br}_{\sxo_{\infty}}})=H^1(C, \ca /\ca_0)\otimes_kH^0(S, \cn_S).$$

On the other hand, the presheaf $S\mapsto \Hom\nolimits_{form. sch.} (S, {\bf \widehat{Br}_{\sxo_{\infty}}})$ is a sheaf in the Zariski topology, since it follows from \cite[1.10.4.6]{EGAI} and the sheaf properties of $\co_U$.

Now let's prove that (\ref{splitex}) is exact.
This sequence is a part of the long exact cohomology sequence that comes from the short exact sequence
$$
0\lto \ca_{S,0}^*\lto \Im_S \lto (\ca /\ca_0)\boxtimes_k \cn_S\lto 0.
$$

Let's show that sequence (\ref{splitex}) is left exact due to our
assumption formula~(\ref{assumpt}). It suffice  to show that the map
$$
H^0(C\times_kS, \ca\widehat{\boxtimes}_k\cn_S)\lto H^0(C\times_kS, (\ca /\ca_0){\boxtimes}_k\cn_S)
$$
is surjective (see formula~(\ref{Ims})), or that the maps
$$
H^0(C\times_kS, \ca_i\widehat{\boxtimes}_k\cn_S)\lto H^0(C\times_kS, (\ca_i /\ca_0){\boxtimes}_k\cn_S)
$$
are surjective for all $i<0$.

For any $i<0$, $h\ge 0$ we define
$$
K_{i,h}\eqdef \Coker (H^0(C,\ca_i /\ca_h)\lto H^0(C,\ca_i/\ca_0)).
$$
We have the following exact sequences:
\begin{equation}
\label{lc}
0\rightarrow H^0(C,\ca_0/\ca_h)\rightarrow H^0(C,\ca_i/\ca_h)\overset{\phi_h}{\rightarrow} H^0(C,\ca_i/\ca_0)\rightarrow K_{i,h}\rightarrow 0,
\end{equation}
\begin{multline}
0\rightarrow H^0(C,\ca_0/\ca_h)\otimes_kH^0(S,\cn_S)\rightarrow H^0(C,\ca_i/\ca_h)\otimes_kH^0(S,\cn_S)\rightarrow \\
 H^0(C,\ca_i/\ca_0)\otimes_kH^0(S,\cn_S)\rightarrow
K_{i,h}\otimes_kH^0(S,\cn_S)\rightarrow 0.
\end{multline}
Since $C$ is a projective curve, the projective systems
$$\{H^0(C,\ca_0/\ca_h)\}_{h\in \sdn} \; \mbox{,} \qquad
\{H^0(C,\ca_0/\ca_h)\otimes_kH^0(S,\cn_S)\}_{h\in \sdn} \;\mbox{,}$$
$$
\{H^0(C,\ca_i/\ca_h)\}_{h\in \sdn} \;\mbox{,} \qquad
  \{H^0(C,\ca_i/\ca_h)\otimes_kH^0(S,\cn_S)\}_{h\in \sdn}$$
 satisfy the ML-condition. Therefore, by \cite[lemma 1]{Ku} the sequences of projective limits are also exact, i.e.,
 using lemma \ref{lamma}, we obtain the exact sequences
\begin{equation}
\label{lp} 0\rightarrow H^0(C,\ca_0)\rightarrow H^0(C,
\ca_i)\rightarrow H^0(C, \ca_i/\ca_0)\rightarrow \limproj_{h \in
\sdn} K_{i,h}\rightarrow 0,
\end{equation}
\begin{multline}
0\rightarrow H^0(C\times_kS, \ca_0\widehat{\boxtimes}_k\cn_S)\rightarrow H^0(C\times_kS,
\ca_i\widehat{\boxtimes}_k\cn_S)\rightarrow \\
H^0(C\times_kS, (\ca_i /\ca_0){\boxtimes}_k\cn_S)\rightarrow
\limproj_h (K_{i,h}\otimes_kH^0(S,\cn_S))\rightarrow 0.
\end{multline}
So, our assertion will follow from the fact that $K_{i,h}=0$.

First, note that $K_i\eqdef \limproj\limits_{h\in \sdn} K_{i,h}=0$.
Indeed,  by the assumption formula~(\ref{assumpt}) and
sequence~(\ref{lp}) we have that
$$
0=\Coker (H^0(C, \ca )\rightarrow H^0(C, \ca /\ca_0))=\limind_i K_i.
$$
We consider the following exact diagram:
$$
\begin{array}{ccccccccc}
&  &&0 &&0 && &\\
&  &&\downarrow &&\downarrow && &\\
0\rightarrow & H^0(C,\ca_0 )& \rightarrow & H^0(C, \ca_i)&\rightarrow &H^0(C,\ca_i/\ca_0)&\rightarrow &K_i &\rightarrow 0\\
& \parallel &&\downarrow &&\downarrow &&\downarrow &\\
0\rightarrow & H^0(C,\ca_0 )& \rightarrow & H^0(C, \ca_{i-1})&\rightarrow &H^0(C,\ca_{i-1}/\ca_0)&\rightarrow &K_{i-1} &\rightarrow 0\\
&  &&\downarrow &&\downarrow && &\\
&  &&H^0(C,\ca_{i-1}/\ca_i) &=&H^0(C,\ca_{i-1}/\ca_i) && &
\end{array}
$$
The diagram chase shows that the map $K_{i}\rightarrow K_{i-1}$ is
injective. Therefore, we must have  $K_i=0$ for any $i<0$.

Now, if $K_{i,h}\neq 0$ for some $h>0$, then this would mean that
$K_i\neq 0$. For, we have $\phi_{h+1}(H^0(C,\ca_i/\ca_{h+1}))\subset
\phi_{h}(H^0(C,\ca_i/\ca_{h}))$ for any $h$ (see
sequence~(\ref{lc})). So, if $K_{i,h}\neq 0$, then
$\phi_{h}(H^0(C,\ca_i/\ca_{h}))\neq H^0(C, \ca_i/\ca_0)$, and a
preimage in $H^0(C, \ca_i/\ca_0)$ of any nonzero element from
$K_{i,h}$ gives a nonzero element in $K_i$. Therefore, $K_{i,h}=0$.

Let's show that (\ref{splitex}) is right exact. This follows from the fact that the map
$$H^2(C\times_kS, \ca_{S,0}^*)\rightarrow H^2(C\times_kS, \Im_S)$$ is an isomorphism.
Indeed, this map is a part of the following diagram, which is exact by lemma \ref{lamma},
corollary \ref{corol1} of this lemma and definitions:
$$
\begin{array}{ccccc}
0\rightarrow &H^2(C\times_kS, \ca_{S,0}^*) & \rightarrow & H^2(C\times_kS, \ca_{S_{red}, 0}^*) & \rightarrow 0\\
&\downarrow &&\parallel &\\
0\rightarrow &H^2(C\times_kS, \Im_S) & \rightarrow & H^2(C\times_kS, \Im_{S_{red}})& \rightarrow 0
\end{array}
$$
(We used here exact sequences:
$$
0 \lto \ca_0 \widehat{\boxtimes}_k \cn_S \lto \ca_{S, 0}^* \lto
\ca_{S_{red}, 0}^* \lto 0 \mbox{,}
$$
$$
0 \lto \ca \widehat{\boxtimes}_k \cn_S \lto \Im_S \lto \Im_{S_{red}}
\lto 0 \mbox{,}
$$
and $\ca_{S_{red}, 0}^* = \Im_{S_{red}}$.)

Now we show that (\ref{splitex}) splits and there is a splitting,
which is  functorial with $S$. We consider the following diagram,
which is exact by lemma~\ref{lamma}, corollary~\ref{corol1} of this
lemma, definitions and our assumptions:
$$
\begin{array}{cccccc}
&0&&0&&\\
&\downarrow &&\downarrow &&\\
 & H^1(C\times_kS,\ca_0\widehat{\boxtimes}_k\cn_S)&\rightarrow &H^1(C\times_kS, \ca_{S,0}^*) & \twoheadrightarrow & H^1(C\times_kS, \ca_{S_{red},0}^*) \\
&\downarrow &&\downarrow&&\parallel \\
 & H^1(C\times_kS,\ca \widehat{\boxtimes}_k\cn_S)&\rightarrow &H^1(C\times_kS, \Im_S) &  \twoheadrightarrow & H^1(C\times_kS, \Im_{S_{red}}) \\
&\downarrow &&\downarrow &&\\
 & H^1(C\times_kS, (\ca /\ca_0)\boxtimes_k\cn_S)& = & H^1(C, \ca /\ca_0)\otimes_kH^0(S, \cn_S) &&\\
&\downarrow &&\downarrow &&\\
&0&&0&&
\end{array}
$$
A splitting of the left vertical exact sequence is given by system
of compatible $k$-linear sections of surjective maps $H^1 (C , \ca_i
/\ca_h) \to H^1 (C , \ca_i /\ca_0)$, $h >0$, $i < 0$ and tensor
multiplication (over $k$) of these sections with identity map on
$H^0(S, \cn_S)$.
  It gives the functorial  with $S$ splitting of
sequence~(\ref{splitex}). (Compare with the proof of
proposition~\ref{prop4}.) The theorem is proved.
\begin{flushright}
$\square$
\end{flushright}

\begin{rem}
\label{ochevidnoe}
The first maps in the rows  of  the last diagram  are embeddings. To show this it suffice to prove that the map
$$H^0(C\times_kS, \ca_{S,0}^*)  \rightarrow  H^0(C\times_kS, \ca_{S_{red},0}^*)$$ is surjective.

If $a\in H^0(C\times_kS, \ca_{S,0})$ is an invertible element, then
its image in $H^0(C\times_kS, \co_{C \times_kS})$ must be also
invertible. Now using corollary~\ref{corola1} of lemma \ref{lamma} and
series for $exp(z)$ and $log(1+z)$, because we assumed $\chara k =0$, we can reduce
the proof to the following fact: the map
$$
H^0(C\times_kS, \ca_{S,1})\lto H^0(C\times_kS, \ca_{S_{red},1})
$$
is surjective. The last fact follows from the following
observations:

 1) using the same arguments as in the proof of lemma~\ref{lamma}, we have
$$
H^0(C\times_kS, \ca_{S,1})\simeq \limproj_{j>1}(H^0(C,\ca_1/\ca_j)\otimes_kH^0(S,\co_S)),\mbox{\quad }
$$
$$
H^0(C\times_kS, \ca_{S_{red},1})\simeq \limproj_{j>1}(H^0(C,\ca_1/\ca_j)\otimes_kH^0(S_{red},\co_{S_{red}}));
$$

2) we have short exact sequences
\begin{multline} \label{mult}
0\rightarrow H^0(C,\ca_1/\ca_j)\otimes_kH^0(S,\cn_S)\rightarrow H^0(C,\ca_1/\ca_j)\otimes_kH^0(S,\co_S) \rightarrow \\
H^0(C,\ca_1/\ca_j)\otimes_kH^0(S_{red},\co_{S_{red}})\rightarrow 0
\end{multline}
for all $j>1$, and the projective system
$\{H^0(C,\ca_1/\ca_j)\otimes_kH^0(S,\cn_S)\}_{j>1}$ satisfies the
ML-condition. Therefore, passing to pojective limit with
respect to $j$ in sequence~(\ref{mult}) we obtain again the short
exact sequence.
\end{rem}

\vspace{0.5cm}

Now we compare the sheaves ${\bf\widetilde{Pic}}_{\sxo_{\infty}}$
and ${\bf Pic}_{\sxo_{\infty}}$. Let ${\bf Pic}^0_{X_{\infty}}$ be
the connected component of zero in the group scheme ${\bf
Pic}_{X_{\infty}}$, which is known to be a closed irreducible
subgroup with $${\bf Pic}^0_{X_{\infty}}(k)=\limproj_{i \ge 0}
Pic^0(X_i)$$ (here $X_i=(C, \ca_0/\ca_{i+1})$ is a scheme), see
\cite[prop.5]{Ku}. Besides, we have the following exact sequence of
sheaves which follows from explicit description of the group scheme
${\bf Pic}_{X_{\infty}}$ in proposition~\ref{reprx}:
$$
0 \lto {\bf Pic}^0_{X_{\infty}} \lto {\bf Pic}_{X_{\infty}} \lto \dz
\lto 0 \mbox{.}
$$

\begin{theo} \label{theolast}
Let $C$ and a ribbon  $\xo_{\infty}$ be as in theorem
\ref{represent}. Assume additionally that $C$ is smooth, and assume
that $\xo_{\infty} $ comes from a smooth projective surface $X$ and
the curve $C\subset X$ with $(C\cdot C)\neq 0$. Then the following
sequence of  sheaves is exact:
$$
0\lto \dz \lto {\bf \widetilde{Pic}}_{\sxo_{\infty}}\lto {\bf
Pic_{\sxo_{\infty}}}\lto 0 \mbox{,}
$$
and ${\bf Pic_{\sxo_{\infty}}}$ is a formal group scheme, which is
non-canonically isomorphic to  $$(\coprod\limits_{i=1}^{|(C\cdot
C)|}{\bf Pic}^0_{X_{\infty}})\times_k{\bf
\widehat{Br}_{\sxo_{\infty}}} \mbox{.}$$
\end{theo}

\Proof. By our assumption, the ribbon satisfies the condition $(**)$. So, by
definition of $\Im_S$ (\ref{Ims}) (see also proposition \ref{orderfn}),
we have the exact sequence of Zariski sheaves for each $S$:
$$
0\lto \Im_S \lto \ca_S^*\lto \dz\lto 0.
$$
Then we have the following exact sequence of Zariski presheaves:
\begin{equation}
\label{presh}
H^1(C\times_kS, \Im_S)\lto H^1(C\times_kS, \ca_S^*)\lto H^1(C\times_kS, \dz ).
\end{equation}
By theorem \ref{zar}, the Zariski sheaf, associated to the presheaf $S\mapsto H^1(C\times_k S, \dz )$, is zero.
(By \cite[vol.I, ch.III, \S 15, th.40, cor.1]{ZS}, the fibres of the morphism $C\times_k S \to S$ are irreducible. Therefore we could apply
theorem~\ref{zar}.) Let's show that the kernel of the first map is  $H^0(C\times_kS, \dz )$. It's enough to prove
that $H^0(C\times_kS, \ca_{S,0}^*)\simeq H^0(C\times_kS, \ca_S^*)$ for a reduced and connected scheme $S$. Indeed, if it is true, then this isomorphism holds for any scheme $S$, because for any affine Noetherian scheme $S$ we then have the following diagram, which is exact by definitions and remark \ref{ochevidnoe}:
$$
\begin{array}{ccccccc}
&&&0&&&\\
&&&\downarrow&&&\\
0\rightarrow &H^0(C\times_kS, \ca\widehat{\boxtimes}_k\cn_S)&\rightarrow &H^0(C\times_kS, \Im_S)&\rightarrow &H^0(C\times_kS, \ca_{S_{red},0}^*)&\rightarrow 1\\
&\parallel &&\downarrow &&\parallel &\\
0\rightarrow &H^0(C\times_kS, \ca\widehat{\boxtimes}_k\cn_S)&\rightarrow &
H^0(C\times_kS, \ca_S^*)&\rightarrow &H^0(C\times_kS, \ca_{S_{red}}^*)&\\
\end{array}
$$

Now, if $S$ is $\Spec k$, the isomorphism $H^0(C\times_kS, \ca_{S,0}^*)\simeq H^0(C\times_kS, \ca_S^*)$ for a reduced $S$ follows from example 8, \cite{Ku}. Recall that in this case we have the exact sequence
$$
0\lto \dz \overset{\alpha}{\lto} Pic (X_{\infty})\lto Pic(\xo_{\infty}),
$$
where $\alpha (1)=\ca_1$, and $\ca_1$ is not a torsion element in the group $Pic (X_{\infty})$, because its image in $Pic (C)$ has degree equal to $-(C\cdot C)\neq 0$. So, the element $\ca_1\widehat{\boxtimes}_k\co_S$ is not a torsion element in $Pic (X_{\infty ,S})$, and therefore $\dz \rightarrow Pic (X_{\infty ,S})$ is injective and $H^0(C\times_kS, \ca_0^*)\simeq H^0(C\times_kS, \ca^*)$ as the long exact sequence
$$
0\rightarrow H^0(C\times_kS, \ca_{S,0}^*) \rightarrow H^0(C\times_kS, \ca_S^*) \rightarrow \dz \rightarrow Pic (X_{\infty ,S})
$$
shows.

Therefore, the sequence (\ref{presh}) leads to the exact sequence of Zariski sheaves
$$
0\lto \dz \lto {\bf \widetilde{Pic}}_{\sxo_{\infty}}\lto {\bf
Pic}_{\sxo_{\infty}}\lto 0.
$$
Immediately from the construction of the sequences above follows
that the sheaf ${\bf Pic}_{X_{\infty}}/\dz $ is representable by the
scheme $\coprod\limits_{i=1}^{|(C\cdot C)|}{\bf
Pic}^0_{X_{\infty}}$. Therefore, using theorem \ref{represent}, we
obtain
$$
{\bf Pic}_{\sxo_{\infty}}\simeq (\coprod\limits_{i=1}^{|(C\cdot
C)|}{\bf Pic}^0_{X_{\infty}})\times_k{\bf
\widehat{Br}_{\sxo_{\infty}}}.
$$
\begin{flushright}
$\square$
\end{flushright}

\begin{rem}
If $C$ is an ample divisor, then  condition~(\ref{assumpt}) from
theorem \ref{represent} is equivalent to the condition
$H^1(X,\co_X)=0$ (see section~\ref{piccoh}).
\end{rem}

\begin{rem}
\label{last}
It would be interesting to obtain the analogoues of
theorems~\ref{represent} and~\ref{theolast} if condition~(\ref{assumpt})
is not satisfied.
It seems that we have to consider  the fpqc sheaves
instead of Zariski sheaves associated with the presheaves
$
S \mapsto \widetilde{Pic}_{\sxo_{\infty}}(S)$
and
 $
S \mapsto {Pic}_{\sxo_{\infty}}(S)$ correspondingly
   to obtain the representability of these sheaves.
\end{rem}

\noindent H. Kurke,  Humboldt University of Berlin, department of
mathematics, faculty of mathematics and natural sciences II, Unter
den Linden 6, D-10099, Berlin, Germany \\ \noindent\ e-mail:
$kurke@mathematik.hu-berlin.de$

\vspace{0.5cm}

\noindent D. Osipov,  Steklov Mathematical Institute, algebra
department, Gubkina str. 8, Moscow, 119991, Russia \\ \noindent
e-mail:
 ${d}_{-} osipov@mi.ras.ru$

\vspace{0.5cm}

\noindent A. Zheglov,  Lomonosov Moscow State  University, faculty
of mechanics and mathematics, department of differential geometry
and applications, Leninskie gory, GSP, Moscow, \nopagebreak 119899,
Russia
\\ \noindent e-mail
 $azheglov@mech.math.msu.su$

\end{document}